\documentclass[11pt]{article}

\usepackage{amstext,amssymb,amsmath,amsbsy,bm}
\usepackage{hyperref}
\usepackage{amscd}
\usepackage{amsfonts}
\usepackage{indentfirst}
\usepackage{verbatim}
\usepackage{amsmath}
\usepackage{amsthm}
\usepackage{enumerate}
\usepackage{graphicx}
\usepackage{color}
\usepackage[OT1]{fontenc}
\usepackage[latin1]{inputenc}
\usepackage[english]{babel}
\usepackage{amssymb}
\usepackage{subfig}
\usepackage{algorithm}
\usepackage{algpseudocode}
\usepackage[shortlabels]{enumitem}

\newcommand{\R}{\mathbb{R}}
\newcommand{\ds}{\displaystyle}

\newcommand{\x}{{\bf x}}

\newcommand{\bv}{{\bf v}}
\newcommand{\bh}{{\bf h}}

\newcommand{\ii}{{\rm i}}
\newcommand{\ik}{\ii k}

\newtheorem{Theorem}{Theorem}[section]
\newtheorem{Lemma}{Lemma}[section]

\newtheorem{Remark}{Remark}[section]

\newtheorem*{Assumption*}{Assumption}

\newtheorem{Problem}{Problem}[section]
\newtheorem*{Problem*}{Problem}
\numberwithin{equation}{section}

\usepackage{matlab-prettifier}
\usepackage{mathtools}
\DeclarePairedDelimiter\ceil{\lceil}{\rceil}

\begin{document}

\title{A global approach for the inverse scattering problem using a Carleman contraction map}

\author{
Phuong M. Nguyen\thanks{
Department of Mathematics and Statistics, University of North Carolina at
Charlotte, Charlotte, NC, 28223, USA, \texttt{pnguye45@charlotte.edu}}
\and
Loc H. Nguyen\thanks{Department of Mathematics and Statistics, University of North Carolina at
Charlotte, Charlotte, NC, 28223, USA, \texttt{loc.nguyen@charlotte.edu}.}
\and
Huong T. T. Vu\thanks{Department of Information Technology, University of Finance-Marketing, Ho Chi Minh City, Vietnam, \texttt{vtthuong@ufm.edu.vn}.} 
}

\date{}
\maketitle
\begin{abstract}
 This paper addresses the inverse scattering problem in a domain $\Omega$. The input data, measured outside $\Omega$,  involve the waves generated by the interaction of plane waves with various directions and unknown scatterers fully occluded inside $\Omega$. 
 The output of this problem is the spatially dielectric constant of these scatterers.
Our approach to solving this problem consists of two primary stages. Initially, we eliminate the unknown dielectric constant from the governing equation, resulting in a system of partial differential equations. Subsequently, we develop the Carleman contraction mapping method to effectively tackle this system.
It is noteworthy to highlight this method's robustness. It does not request a precise initial guess of the true solution, and its computational cost is not expensive.
Some numerical examples are presented.
\end{abstract}

\noindent{\it Keywords: inverse scattering problem; Carleman contraction mapping method; globally convergence; low computational cost.}

\noindent{\it AMS subject classification: 	35R25, 35N10, 35R30, 78A46

}

\section{Introduction} \label{sec intr}

\begin{figure}[h!]
\begin{center}
	\includegraphics[width=.4\textwidth]{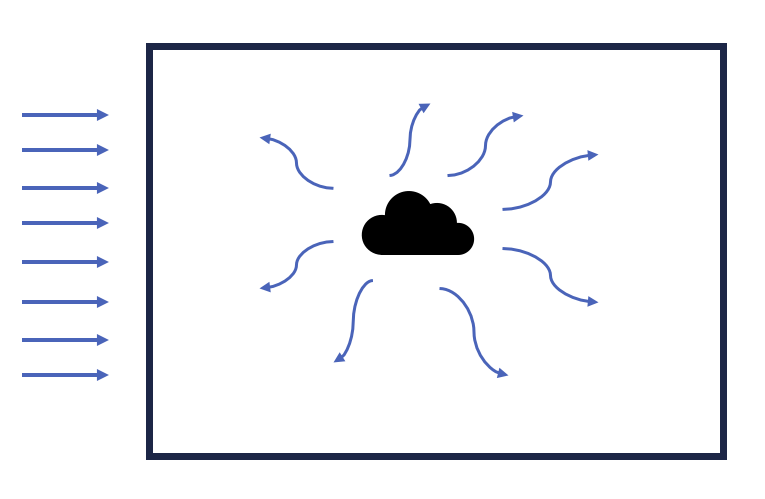}
	\put(-100,57){\tiny A unknown scatterer}
	\put(-70,34){\tiny Scattering wave}
	\put(-110,104){\tiny Scattering wave}
	\put(-200,85){\tiny \rotatebox{270}{Incident wave}}
	\put(-154,15){\tiny $\Omega$}
	\caption{\label{figdia} A visualization of the inverse scattering process: An incident wave hits a scatterer in an inaccessible domain $\Omega$. The interaction between the wave and the scatterer produces a scattered wave, which contains important information about the scatterer (shape, location, and some physical properties). The primary aim of the inverse scattering problem is to use the measurements of the scattering wave outside $\Omega$ to find this information.}
\end{center}
\end{figure}
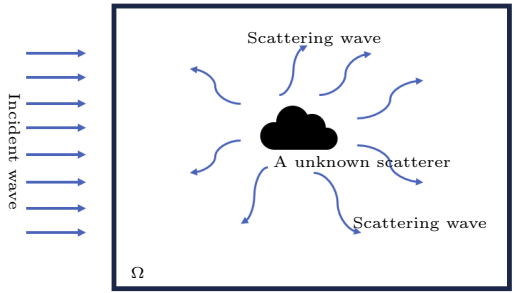

Let $d \geq 2$ be the spatial dimension.
Let $c: \R^d \to \R$ be the spatially distributed dielectric constant of the medium. Let $\Omega$ be a domain of $\R^d$ with a smooth boundary.
Assume that $c$ satisfies
\begin{equation}
	c(\x) =
	\left\{
		\begin{array}{ll}
			\geq 1 &\x \in \Omega,\\
			1 &\x \in \R^d \setminus \Omega.
		\end{array}
	\right.
	\label{1.1}
\end{equation}
Setting $c(\x) = 1$ outside $\Omega$ means that $\R^d \setminus \Omega$ is either vacant or air-filled, whereas $c(\x) \not= 1$ inside $\Omega$ indicates the presence of entities such as obstacles, acting as scatterers. To identify these scatterers, we use optical waves (known as the incident wave) to illuminate $\Omega$. These incident waves propagate through space and scatter in all directions. By recording the scattered waves and using these data to calculate the spatially distributed dielectric constant $c(\x),$ $\x \in \Omega$, we can map the internal geometry of $\Omega$. The knowledge of this function $c$ helps identify the shape, position, and certain optical characteristics of the obstacles we aim to detect.
Typical examples of those scatterers are anti-personnel explosive devices buried under the ground, cancerous tumors inside living tissues, and nanostructures.
Therefore, the success in finding $c$ strongly impacts many other fields such as bio-medical imaging, nondestructive testing, radar, detection and identification of explosive devices buried under the ground, optical physics, seismic exploration, and nanoscience.
 See Figure \ref{figdia} for a 2D diagram of the scattering phenomenon.
In this paper, we use plane waves with a fixed wave number $k$ and various angles $\theta \in [0, 2\pi]$ as the incident waves.
The incident plane wave is given by
\begin{equation}
	u_{\rm inc}(\x, \theta) = e^{\ik \hat \theta \cdot \x},
	\quad 
	\mbox{for all } \x \in \R^d,
\end{equation}
where 
\begin{equation}
	\hat \theta = (\cos \theta, 0, \dots, 0, \sin\theta), \quad \theta \in [0, 2\pi].
	\label{1.3}	
\end{equation} 
The resulting total wave $u = u(\x, \theta)$, $(\x, \theta) \in \R^d \times [0, 2\pi]$, satisfies the Helmholtz equation and the Sommerfeld radiation condition
\begin{equation}
	\left\{
		\begin{array}{ll}
			\Delta u + k^2 c(\x) u = 0 &\x \in \R^d,\\
			\partial_{|\x|} (u - u_{\rm inc}) - \ik (u - u_{\rm inc}) = o(|\x|^{\frac{1 - d}{2}}) &|\x| \to \infty.
		\end{array}
	\right.
	\label{Hel_eqn}
\end{equation}
We are interested in the following inverse scattering problem:
\begin{Problem}[Inverse scattering problem]
	Given the boundary data
	\begin{equation}
		f(\x, \theta) = u(\x, \theta)
		\quad
		\mbox{and}
		\quad
		g(\x, \theta) = \partial_{\nu} u(\x, \theta)
		\label{1.5}
	\end{equation}
	for all $(\x, \theta) \in \partial \Omega \times [0, 2\pi],$ determine the function $c(\x),$ for all $\x \in \Omega.$
	\label{p}
\end{Problem}

It is worth mentioning that the data in \eqref{1.5} are suitably requested. The target function to be determined, $c(\x)$ has $d$ variables $x_1, \dots, x_d$. The domain $\partial \Omega \times [0, 2\pi]$ of the data has a combined dimension of $(d - 1) + 1$.  From a mathematical standpoint, there is a precise alignment between the dimensions of the required data and that of the target function. Consequently, Problem \ref{p} is efficiently designed, requiring no redundant data, and is not over-determined. Our procedure to solve Problem \ref{p} consists of two steps.
 In the first step, we eliminate the unknown coefficient $c$ from \eqref{Hel_eqn}. The result of this step is a system of partial differential equations (PDEs), namely $\mathcal L \bv = 0$ where $\mathcal L$ is a nonlinear differential operator.  
In the second step, we solve the system $\mathcal L \bv = 0$ for the vector-valued function $\bv$.
Once the solution $\bv$ to the system above is found, we can compute the solution to the governing equation.
The solution to the inverse scattering problem follows.
The first step is based on a change of variable.
 The second step is challenging since $\mathcal L$ is highly nonlinear. 
Rather than using the optimization-based method, which requires a good initial guess of the true solution, we will employ the recently developed method in \cite{AbneyLeNguyenPeters, LeCON2023, LeNguyen:jiip2022, Nguyen:AVM2023} to find $\bv$. More precisely, we construct a contraction mapping whose unique fixed point is an approximation of $\bv$. The solution $c$ to Problem \ref{p} follows.

 There are many effective approaches to solve various versions of the inverse scattering problem:
 \begin{itemize}
 \item  the imaging techniques based on sampling and the factorization methods \cite{ AmmariChowZou:sjap2016, AmmariKang:lnim2004,  Tan1:ip2012, Burge2005,   Colto1996, HarrisLiem:SIAM2020, Kirsc1998, LiLiuZou:smms2014, LiLiuWang:jcp2014, Liem:jiip2022, Soumekh:SAR};
 \item the methods based on optimization \cite{Bakushinsii:springer2004, Tan2:ip2012, Chavent:springer2009, Engl:Kluwer1996, Gonch2013, Tihkonov:kapg1995};
 \item the methods based on Born series \cite{Bleistein:ap1984, Chew:vnr1990, Devaney:cup2012, Kirsch:aa2017, Langenberg:1987, Moskow:Ip2008};
 \item the method based on linearization \cite{BaoLi:ip2005, BaoLi:SIAM2005, Bao:ip2015, Chen:ip1997}; 
 \item  the convexification method \cite{VoKlibanovNguyen:IP2020, Khoaelal:IPSE2021, KhoaKlibanovLoc:SIAMImaging2020};
 \item
see \cite{ColtonKress:2013} for a relatively complete list.
\end{itemize}
Although effective, each approach above has drawbacks.
For example, the imaging techniques only provide  locations and shapes of the unknown objects; the methods based on optimization are local in the sense that they require a good initial guess of the true solution.
There are many existing methods to globalize the local algorithms,
one of which is the convexification method, first introduced in \cite{KlibanovIoussoupova:SMA1995}.
The key of the convexification method is to include a Carleman weight function into the mismatch functional.
Using Carleman estimates, we can prove that the mismatch functional with the presence of suitable Carleman weight functions are strictly convex.
Hence, the mismatch functional has a unique minimizer.
We can rigorously prove that this unique minimizer is an approximation of the true solution of the inverse problem as the noise contained in the data tends to $0$ with the Lipchitz rate.
Our research team has developed new versions of the convexification method to solve inverse scattering problem in the frequency domain \cite{VoKlibanovNguyen:IP2020, Khoaelal:IPSE2021,KhoaKlibanovLoc:SIAMImaging2020} 
and in the time domain \cite{KlibanovLeNguyenIPI2022}.
These methods based on convexification were verified by both simulated and experimental data. 
In this paper,  as mentioned in the previous paragraph, we will employ the contraction principle rather than convexifying the mismatch functional whose unique minimizer yields the desired solution. Given that the recursive iteration procedure to compute that fixed point converges at a linear rate, by proposing this method, we will improve the reconstruction quality and to reduce the computational time.

The paper is organized as follows. In Section \ref{sec2}, we derive an approximation model, which is a system of PDEs. In Section \ref{sec3}, we define a contraction mapping and propose an iterative procedure to solve that system. 
In Section \ref{sec4}, we present some numerical examples and some details in implementation.
Section \ref{sec5} is for some discussions and concluding remarks.

\section{An approximate model}\label{sec2}

In this section, we outline the derivation of a system of PDEs, the solution of which provides the knowledge of $c$.
While this step was introduced and used previously in \cite{VoKlibanovNguyen:IP2020, Khoaelal:IPSE2021,KhoaKlibanovLoc:SIAMImaging2020}, we recall them here for completeness and for the reader's convenience. 
	Introduce the change of variable
	\begin{equation}
		v(\x, \theta) = \frac{1}{k^2}\log \frac{u(\x, \theta)}{u_{\rm inc}(\x, \theta)}
		\quad
		\mbox{for all } (\x, \theta) \in \Omega \times [0, 2\pi].
		\label{change_variable}
	\end{equation}
	Defining the logarithm of the complex function $\frac{u(\x, \theta)}{u_{\rm inc}(\x, \theta)}$, $(\x, \theta) \in \Omega \times [0, 2\pi]$, in \eqref{change_variable} presents a challenge. However, this logarithm is well-defined under some certain circumstances. For instance, according to \cite[Theorem 3]{KlibanovLiemLoc:SIAMjap2019}, the function $\frac{u(\x, \theta)}{u_{\rm inc}(\x, \theta)}$ can be approximated by $A(\x)e^{\ik (\tau(\x) - \x\cdot \hat \theta)} + O(1/k)$ when $k$ is sufficiently large, where $A$ and $\tau$ are some positive functions. It is important to note that \cite[Theorem 3]{KlibanovLiemLoc:SIAMjap2019} was proved for a 3D plane wave in the form of $e^{\ik x_3}$. However, with a variable change (if necessary), this theorem can be applied to the incident wave in the form $e^{\ik \x \cdot \hat \theta}$ with any direction as used in this paper. More significantly, according to \cite[Theorem 3.3]{Klibanov:jiip2023}, the function $A$ exceeds a positive constant $a_0$. Therefore, the function $v$ in \eqref{change_variable} is well-defined and can be approximated by $\frac{1}{k^2}[\ln A(\x) + \ik (\tau(\x) - \x \cdot \hat \theta)]$. In the numerical implementation, we defined the logarithm of the above function straightforwardly, encountering no difficulties.

	By a straightforward algebra, we have
	\begin{align}
		\nabla v(\x, \theta) 
		&= \frac{u_{\rm inc}(\x, \theta)}{k^2 u(\x, \theta)} \nabla \frac{u(\x, \theta)}{u_{\rm inc}(\x, \theta)} 
		= \frac{u_{\rm inc}(\x, \theta)}{k^2 u(\x, \theta)} 
		\frac{u_{\rm inc}(\x, \theta) \nabla u(\x, \theta) - u(\x, \theta) \nabla u_{\rm inc}(\x, \theta)}{u_{\rm inc}^2(\x, \theta)} \nonumber\\
		&= \frac{1}{k^2}\Big[
			\frac{\nabla u(\x, \theta)}{u(\x, \theta)}
			- \frac{\nabla u_{\rm inc}(\x, \theta)}{u_{\rm inc}(\x, \theta)}
		\Big]
		= \frac{1}{k^2} \Big[
			\frac{\nabla u(\x, \theta)}{u(\x, \theta)}
			- \ik \hat \theta
		\Big]
		\label{2.1}
	\end{align}
	and
	\begin{equation}
		\Delta v(\x, \theta)
		= \frac{1}{k^2} \Big[
			\frac{\Delta u(\x, \theta)}{u(\x, \theta)}
			- \Big(\frac{\nabla u(\x, \theta) }{u(\x, \theta)}\Big)^2
		\Big]
		\label{2.2}
	\end{equation}
	for all $(\x, \theta) \in \Omega \times [0, 2\pi]$.
	It follows by the Helmholtz equation in \eqref{Hel_eqn} that $\frac{\Delta u(\x, \theta)}{u(\x, \theta)} = -k^2c(\x)$.
	On the other hand, by \eqref{2.1}, $\frac{\nabla u(\x, \theta)}{u(\x, \theta)} = k^2 \nabla v(\x, \theta) + \ik \hat\theta$.
	Combining these identities and \eqref{2.2} gives
	\begin{align*}
		\Delta v(\x, \theta) 
		&= 
		\frac{1}{k^2}\Big[
			-k^2 c(\x) - (k^2\nabla v(\x, \theta) + \ik \hat\theta)^2
		\Big]\\
		&=-c(\x) - k^2(\nabla v(\x, \theta))^2 - 2\ik \nabla v(\x, \theta) \cdot \hat \theta
		+ 1
	\end{align*}
	for all $(\x, \theta) \in \Omega \times [0, 2\pi]$.
	Therefore, 
	the function $v$ is governed by the equation
	\begin{equation}
		\Delta v(\x, \theta) + 2\ik \nabla v(\x, \theta) \cdot \hat \theta + k^2(\nabla v(\x, \theta))^2  = -c(\x) + 1		
		\label{2.3}
	\end{equation}
	for all $(\x, \theta) \in \Omega \times [0, 2\pi]$.
	Differentiate equation \eqref{2.3} with respect to $\theta$ and note (from \eqref{1.3}) that 
	\begin{equation}
		\frac{d \hat \theta}{d\theta} = \hat \theta' =  (-\sin \theta, 0, \dots, 0, \cos\theta) 
				\quad
		\mbox{for all } \theta \in [0, 2\pi].
	\end{equation}
	We have
	\begin{equation}
		\Delta v_\theta(\x, \theta) + 2\ik \nabla v_\theta(\x, \theta) \cdot \hat \theta
		+
		2\ik \nabla v(\x, \theta) \cdot \hat \theta' + 2k^2\nabla v(\x, \theta) \cdot \nabla v_\theta(\x, \theta)
		=
		0
		\label{2.5}
	\end{equation}
	for $(\x, \theta) \in \Omega \times [0, 2\pi]$.
	
Let $\{\Psi\}_{n \geq 1}$ be the polynomial-exponential orthonormal basis of $L^2(0, 2\pi)$.
This basis was constructed in \cite{Klibanov:jiip2017} by using the Gram-Schmidt process on the complete set $\{\phi_n(\theta) = \theta^{n-1}e^\theta\}_{n \geq 1}$ of $L^2(0, 2\pi)$. See Remark \ref{rem 21} for the reason why we choose this basis among many other orthonormal bases of $L^2(0, 2\pi)$.
For each $\x \in \Omega$, we write the Fourier expansion of $v(\x, \theta)$ with respect to the basis $\{\Psi\}_{n \geq 1}$
\begin{equation}
	v(\x, \theta) = \sum_{n= 1}^\infty v_n(\x) \Psi_n(\theta),
	\quad
	\mbox{where } v_n(\x) = \int_{0}^{2\pi}v(\x, \theta) \Psi_n(\theta) d\theta.
	\label{2.6}
\end{equation}
Motivated by the Galerkin approximation and since we pay attention on computation, we approximate $v(\x, \theta)$ by truncating the series in \eqref{2.6}
\begin{equation}
	v(\x, \theta) = \sum_{n = 1}^N v_n(\x) \Psi_n(\theta),
	\label{2.7}
\end{equation}
for $(\x, \theta) \in \Omega \times [0, 2\pi]$,
for some cut-off number $N$. We will discuss a data-driven approach to choose $N$ later in Subsection \ref{sub42}.
Plugging \eqref{2.7} into \eqref{2.5}, we have
\begin{multline}
		\sum_{n = 1}^N\Delta  v_n(\x) \Psi_n'(\theta) + 2\ik  \sum_{n = 1}^N \nabla  v_n(\x) \Psi_n'(\theta) \cdot \hat \theta
		+
		2\ik \sum_{n = 1}^N \nabla  v_n(\x) \Psi_n(\theta) \cdot \hat \theta' 
		\\
		+ 2k^2\sum_{n = 1}^N \sum_{l = 1}^N \nabla v_n(\x)  \cdot  \nabla v_l(\x) \Psi_n(\theta)\Psi_l'(\theta)
		=
		0
		\label{2.8}
	\end{multline}
	for $(\x, \theta) \in \Omega \times [0, 2\pi]$.
	For each $m \in \{1, \dots, N\}$, we multiply $\Psi_m(\theta)$ to both sides of \eqref{2.8}.
	Then, integrating the resulting equation, we obtain
	\begin{equation}
		\sum_{n = 1}^Ns_{mn}\Delta  v_n(\x) 
		+  \sum_{n = 1}^N \nabla  v_n(\x)  \cdot B_{mn}		 
		+ \sum_{n = 1}^N \sum_{l = 1}^N a_{mnl}\nabla v_n(\x)  \cdot  \nabla v_l(\x) 
		=
		0
		\label{2.9}
	\end{equation}
	for all $\x \in \Omega$
	where
	\begin{align*}
		s_{mn} &= \int_{0}^{2\pi}\Psi_n'(\theta)\Psi_m(\theta) d\theta,\\
		B_{mn} &= 2\ik\Big[\int_{0}^{2\pi}   \hat \theta \Psi_n'(\theta) \Psi_m(\theta)d\theta
		+ \int_{0}^{2\pi}  \hat \theta' \Psi_n(\theta)  \Psi_m(\theta) d\theta\Big],
		\\
		a_{mnl} &= 2k^2\int_{0}^{2\pi} \Psi_n(\theta)\Psi_l'(\theta) \Psi_m(\theta)d\theta.
	\end{align*}
	Due to the presence of $s_{mn}$, the system \eqref{2.9} might not be elliptic. However, it is equivalent to an elliptic system.
	Recall that the matrix $S$ is invertible and denote by $(\widetilde s_{mn})_{m, n = 1}^N$ the matrix $S^{-1}$.
	Define $\widetilde  B_{mn} = \sum_{i = 1}^N\widetilde s_{mi} B_{in}$ and $\widetilde a_{mnl} = \sum_{i = 1}^N\widetilde s_{mi} a_{inl}$ for $m, n, l \in \{1, \dots, N\}.$
	It follows from \eqref{2.9} that
	\begin{equation}
		\Delta  v_m(\x) 
		+  \sum_{n = 1}^N \nabla  v_n(\x)  \cdot \widetilde  B_{mn}		 
		+ \sum_{n = 1}^N \sum_{l = 1}^N \widetilde  a_{mnl}\nabla v_n(\x)  \cdot  \nabla v_l(\x) 
		=
		0
		\label{2.99}
	\end{equation}
	for $\x\in \Omega$.	
	Coupling \eqref{2.99} for $m \in \{1, \dots, N\}$, we obtain a nonlinear system for the vector $\bv = 
	[
		\begin{array}{lll}
			v_1 &\dots &v_N
		\end{array}
	]^{\rm T}
	$ where the superscript $^{\rm T}$ indicates the transpose.

\begin{Remark}
The system formed by \eqref{2.9} with $m \in \{1, 2, \dots, N\}$ plays the key role in our algorithm. It is our approximation model for Problem \ref{p}.
The choice of $\{\Psi\}_{n \geq 1}$ is important in establishing \eqref{2.9}.
For any $N \in \mathbb{N}$, let $S$ be the $N \times N$ matrix $S$ whose $mn-$entry is $s_{mn} = \int_{0}^{2\pi}\Psi_n'(\theta)\Psi_m(\theta) d\theta$ as above. It was proved in \cite{Klibanov:jiip2017} that:
\begin{enumerate}
	\item the matrix $S$ is invertible;
	\item $\Psi'_n$ is not identically zero in $[0, 2\pi],$ for all $n \geq 1$.
\end{enumerate}

It is obvious that the invertibility of $S$ is important. Without the knowledge of the matrix $S^{-1}$, defining coefficients for \eqref{2.9} is impossible.
On the other hand,
the second property mentioned above is also holds significant. Specifically, if $\Psi'_n = 0$ for some $n$, crucial information regarding $\Delta v_n$ in the first term of \eqref{2.8} would be lost.
To illustrate, consider replacing the polynomial-exponential basis with widely used alternatives such as the Legendre polynomial or the trigonometric basis, denoted by $\{\phi_n\}_{n \geq 1}$. In this scenario, $\phi_1$ is a constant, resulting in $\phi_1' \equiv 0$. Consequently, the principal part $\Delta v_1$ would not be present in \eqref{2.8} and, subsequently, in \eqref{2.9}.
As a consequence, the error in computing $v_1$ might be substantial, potentially leading to inaccuracies in the recovered $v(\x, \theta)$ using \eqref{2.7}. This is primarily because the first term holds significant importance in this truncated series.
\label{rem 21}
\end{Remark}

	We next compute the boundary information for the vector $V$. 
	For each $m \in \{1, \dots, N\},$ due to \eqref{2.6},
	\begin{equation}
		v_m(\x) = \int_{0}^{2\pi} v(\x, \theta)\Psi_m(\theta) d\theta
		= \frac{1}{k^2}\int_{0}^{2\pi} \log \frac{u(\x, \theta)}{u_{\rm inc}(\x, \theta)}\Psi_m(\theta) d\theta = \frac{1}{k^2}\int_{0}^{2\pi} \log \frac{f(\x, \theta)}{u_{\rm inc}(\x, \theta)}\Psi_m(\theta) d\theta
		\label{2.11}
	\end{equation}
	for all $\x \in \partial \Omega.$
	Here, we have used \eqref{change_variable} and the knowledge of the measurement $f(\x, \theta) = u(\x, \theta)$ on $\partial \Omega \times [0, 2\pi]$ as in \eqref{1.5}.
	Similarly, using \eqref{2.1}, we can compute the normal derivative of $v_m$, $m \in \{1, \dots, N\}$, as
	\begin{equation}
		\partial_{\nu} v_m(\x) = \frac{1}{k^2} \int_{0}^{2\pi} \Psi_m(\theta) \Big[
			\frac{\nabla u(\x, \theta) \cdot \nu }{u(\x, \theta)}
			- \ik \hat \theta \cdot \nu 
		\Big]  d\theta
		= 
		 \frac{1}{k^2} \int_{0}^{2\pi} \Psi_m(\theta) \Big[
			\frac{g(\x, \theta)}{f(\x, \theta)}
			- \ik \hat \theta \cdot \nu 
		\Big]  d\theta
	\end{equation}
	for all 
	$\x \in \partial \Omega.$
Defining 
\begin{align}
	F_m(\x) & = \frac{1}{k^2}\int_{0}^{2\pi} \log \frac{f(\x, \theta)}{u_{\rm inc}(\x, \theta)}\Psi_m(\theta) d\theta, \label{Fm}\\
	G_m(\x) &= \frac{1}{k^2} \int_{0}^{2\pi} \Psi_m(\theta) \Big[
			\frac{g(\x, \theta)}{f(\x, \theta)}
			- \ik \hat \theta \cdot \nu 
		\Big]  d\theta \label{Gm}
\end{align}
	for $\x \in \partial \Omega$, $m \in \{
		1, \dots, N
	\}$, we obtain a system for $\bv = 
	[
		\begin{array}{lll}
			v_1 &\dots &v_N
		\end{array}
	]^{\rm T}
	$
\begin{equation}
	\left\{
		\begin{array}{ll}
		\ds \Delta  v_m(\x) 
		+  \sum_{n = 1}^N \nabla  v_n(\x)  \cdot \widetilde B_{mn}		 
		+ \sum_{n = 1}^N \sum_{l = 1}^N \widetilde  a_{mnl}\nabla v_n(\x)  \cdot  \nabla v_l(\x) 
		=
		0
		&\x \in \Omega,\\
		v_m(\x) = F_m(\x) &\x \in \partial \Omega,\\
		\partial_{\nu} v_m(\x) = G_m(\x) &\x \in \partial \Omega
		\end{array}
	\right.
	\label{2.13}
\end{equation}
	for $m \in \{1, \dots, N\}.$

The derivation of \eqref{2.13} suggests a numerical method to solve Problem \ref{p}.
Initially, we compute the solution to the nonlinear problem \eqref{2.13}. Next, we use \eqref{2.7} to evaluate $v(\x, \theta)$ for $(\x, \theta) \in \Omega \times [0, 2\pi].$
Lastly, the desired solution $c(\x)$ can be computed using \eqref{2.3}.
Tackling \eqref{2.13} poses considerable difficulty due to its nonlinear nature, prompting the application of the Picard iteration method. 
Let ${\bf v}^{(0)}$ be an ``initial guess" and assume by induction that the vector ${\bf v}^{(p)} = [
		\begin{array}{lll}
			v_1^{(p)} &\dots &v_N^{(p)}
		\end{array}
	]^{\rm T} $, $p \geq 0$, is given.
We compute ${\bf v}^{(p+1)} = [
		\begin{array}{lll}
			v_1^{(p+1)} &\dots &v_N^{(p+1)}
		\end{array}
	]^{\rm T}$ by solving
\begin{equation}
	\left\{
		\begin{array}{ll}
			\ds \Delta  v_m^{(p + 1)}(\x) 
		+  \sum_{n = 1}^N \nabla  v_n^{(p + 1)}(\x)  \cdot \widetilde B_{mn}		
		\\ 
		\hspace{5cm} \ds + \sum_{n = 1}^N \sum_{l = 1}^N \widetilde a_{mnl}\nabla v_n^{(p + 1)}(\x)  \cdot  \nabla v_l^{(p)}(\x) 
		=
		0
		&\x \in \Omega,\\
		v_m^{(p + 1)}(\x) = F_m(\x) &\x \in \partial \Omega,\\
		\partial_{\nu} v_m^{(p + 1)}(\x) = G_m(\x) &\x \in \partial \Omega
		\end{array}
	\right.
	\label{2.14}
\end{equation}
Given that \eqref{2.14} is linear with respect to ${\bf v}^{(p+1)}$, it can be solved through the quasi-reversibility method or the least-squares optimization. A critical consideration is the convergence of the sequence $\{{\bf v}^{(p)}\}_{p \geq 0}$ to the actual solution ${\bf v}^*$ of \eqref{2.13}, particularly when the initial guess ${\bf v}^{(0)}$ far away from ${\bf v}^*$.
To guarantee the convergence, we will include a Carleman weight function in the least-squares optimization above. Details about combining a Carleman weight function and the least squares optimization will be provided in the next section.

\begin{Remark}
	The Picard iterative scheme above is different from the ones suggested in \cite{LeCON2023, LeNguyen:jiip2022, Nguyen:AVM2023}. 
	In fact, if we applied the scheme in those papers, problem \eqref{2.14} would be read as
	\begin{equation}
	\left\{
		\begin{array}{ll}
			\ds\Delta  v_m^{(p + 1)}(\x) 
		+  \sum_{n = 1}^N \nabla  v_n^{(p )}(\x)  \cdot \widetilde  B_{mn}		
		 + \sum_{n = 1}^N \sum_{l = 1}^N \widetilde  a_{mnl}\nabla v_n^{(p)}(\x)  \cdot  \nabla v_l^{(p)}(\x) 
		=
		0
		&\x \in \Omega,\\
		v_m^{(p + 1)}(\x) = F_m(\x) &\x \in \partial \Omega,\\
		\partial_{\nu} v_m^{(p + 1)}(\x) = G_m(\x) &\x \in \partial \Omega.
		\end{array}
	\right.
	\label{2.1414}
\end{equation}
In contrast, we solve $\bv^{(p+1)}$ by solving \eqref{2.14} rather than \eqref{2.1414} because the nonlinear operator
\[
	F_{\rm nonlinear}(\bv) = \sum_{n = 1}^N \sum_{l = 1}^N \widetilde  a_{mnl}\nabla v_n(\x)  \cdot  \nabla v_l(\x),
	\quad
	\bv \in H,
\]
	is not Lipschitz, as requested in \cite{LeCON2023, LeNguyen:jiip2022, Nguyen:AVM2023}.
	Since the Lipschitz condition is not satisfied, we numerically experience the divergence in computation when applying the Picard scheme in those papers.
	Therefore, the change of replacing \eqref{2.1414} by \eqref{2.14} is significant to the success to our numerical scheme. 
	\label{rem 22}
\end{Remark}

\section{The Carleman contraction mapping method for \eqref{2.13}} \label{sec3}

The first version of the Carleman contraction mapping method was introduced in \cite{LeNguyen:jiip2022}.
We refer the reader to the following works \cite{AbneyLeNguyenPeters, LeCON2023, LeNguyenNguyenPark, NguyenNguyenTruong:camwa2022, Nguyen:AVM2023} for following up works and their application in several different inverse problems for PDEs.
The main idea of this method is to combine the well-known Picard iteration, defined in the proof of the contraction mapping principle, and a suitable Carleman estimate.
For convenience, we recall here a Carleman estimate, which was established in \cite{Nguyen:AVM2023}, which is important for the convergence of our iterative method.

\begin{Lemma}[Carleman estimate]
Fix a point $\x_0 \in \R^d \setminus \Omega$. 
	Define $r(\x) = |\x - \x_0|$ for all $\x \in \Omega.$
There exist  positive constants $\beta_0$ depending only on $\x_0$, $\Omega$, and $d$ such that for all function $v \in C^2(\overline \Omega)$ satisfying 
	\begin{equation}
		v(\x) = \partial_{\nu} v(\x) = 0 \quad \mbox{for all } 
		\x \in \partial \Omega,
		\label{3.1}
	\end{equation}
	the following estimate holds true
	\begin{equation}
		\int_{\Omega} e^{2\lambda r^{-\beta}(\x)}\vert \Delta v \vert^2d\x
		\geq
		 C\lambda  \int_{\Omega}  e^{2\lambda r^{-\beta}(\x)}\vert\nabla v(\x)\vert^2\,d\x
		+ C\lambda^3  \int_{\Omega}   e^{2\lambda r^{-\beta}(\x)}\vert v(\x)\vert^2\,d\x	
		\label{Car est}	
	\end{equation}
	for all $\beta > \beta_0$ and $\lambda \geq \lambda_0$. 
	Here, $\lambda_0 = \lambda_0( \x_0,  \Omega, d, \beta)$ and $C = C( \x_0,  \Omega, d, \beta) > 0$ depend only on the listed parameters.
	\label{lemma carl}
\end{Lemma}

The derivation of Lemma \ref{lemma carl} is directly based on \cite[Lemma 5]{MinhLoc:tams2015}, see \cite[Lemma 2.1]{LeNguyenTran:CAMWA2022} for a comprehensive explanation. Alternatively, another method to obtain \eqref{Car est} involves using a distinct Carleman weight function through the Carleman estimate presented in \cite[Chapter 4, Section 1, Lemma 3]{Lavrentiev:AMS1986}, which applies to general parabolic operators.
It is important to point out to the reader the diverse spectrum of Carleman estimates for different types of differential operators (elliptic, parabolic, and hyperbolic) and their application in the fields of inverse problems and computational mathematics, as discussed in sources such as \cite{BeilinaKlibanovBook, BukhgeimKlibanov:smd1981, KlibanovLiBook, LocNguyen:ip2019}.
Moreover, it is notable that certain Carleman estimates are effective for any functions $v$ meeting the conditions $v\vert_{\partial \Omega} = 0$ and $\partial_{\nu} v\vert_{\Gamma} = 0$, where $\Gamma$ is a subset of $\partial \Omega$. This is exemplified in works like \cite{KlibanovNguyenTran:JCP2022, NguyenLiKlibanov:2019}. These Carleman estimates were especially useful in tackling quasilinear elliptic PDEs with data limited to the portion $\Gamma$ of $\partial \Omega$.

Let $s$ be an integer with $s > \ceil{\frac{d}{2}} + 2$ such that $H^s(\Omega)$ can be continuously embedded into $C^2(\overline \Omega)$.
For the theoretical purpose, we assume that the target function $c$ is sufficiently smooth such that the true solution to the forward problem $u$ is smooth with $\|u(\cdot, \theta)\|_{H^s(\overline \Omega)}$ is finite uniformly in $\theta \in [0, 2\pi].$ We also assume that $|u(\x, \theta)|$ is uniformly bounded from below by a positive constant $u_0$.  Due to the change of variable \eqref{change_variable}, the function $v(\x, \theta)$ is in the class $H^s(\overline \Omega)$ with $H^s$ norm being uniformly bounded for all $\theta \in [0, 2\pi].$ Hence, there is a number $M$, depending on an upper bounded of $\|u(\cdot, \theta)\|_{H^s(\Omega)}$, $\theta \in [0, 2\pi]$ and the set $\{\Psi_n\}_{n = 1}^N$, such that $\|{\bf v}^*\|_{H^s(\overline \Omega)^N} < M$ where ${\bf v}^*$ is the true solution to \eqref{2.13}.
We find the solution to \eqref{2.13} in the set of admissible solutions
	\begin{equation}
		H = 
		\big\{
			\bv = [
		\begin{array}{lll}
			v_1 &\dots &v_N
		\end{array}
	]^{\rm T} \in H^s(\Omega)^N: 
	\|\bv\|_{C^2(\overline \Omega)^N}  \leq M, 
	v_m|_{\partial \Omega} = F_m \mbox{ and } 
	\partial_{\nu} v_m|_{\partial \Omega} = G_m
		\big\}.
		\label{3.3}
	\end{equation}
	Throughout the paper, we assume that $H$ is nonempty and the true solution $v^*$ to \eqref{2.13} belongs to $H$.

As outlined in the paragraph before Remark \ref{rem 22} in Section \ref{sec2}, we construct a sequence $\{{\bf v}^{(p)}\}_{p \geq 0}$ with the hope that it converges to the true solution to \eqref{2.13}. Starting with $\bv^{(0)}$ as a vector-valued function in $H$, we proceed under the induction assumption that $\bv^{(p)}$, $p \geq 0$, is already determined. The subsequent vector $\bv^{(p + 1)}$ is then derived by solving \eqref{2.14}. Given that \eqref{2.14} is linear with respect to $\bv^{(p + 1)}$,  the least squares method (so-called the quasi-reversibility method in solving PDEs) is applicable to compute $\bv^{(p + 1)}$. To guarantee global convergence, it is important to include the Carleman weight function $e^{2\lambda r^{-\beta}(\x)}$ in the corresponding least-squares cost functional. More precisely, we set
\begin{equation}
	\bv^{(p + 1)} =  \underset{\varphi \in H}{\rm argmin} J^{\bv^{(p)}}_{\lambda, \beta, \epsilon}(\varphi)
	\label{3.4}
\end{equation}
where 
\begin{multline}
	J^{\bv^{(p)}}_{\lambda, \beta, \epsilon}(\varphi) = \sum_{m = 1}^N\int_{\Omega} e^{2\lambda r^{-\beta}(\x)}\Big| \Delta  \varphi_m(\x) 
		+  \sum_{n = 1}^N \nabla  \varphi_n(\x)  \cdot \widetilde B_{mn} 
		\\
		+ \sum_{n = 1}^N \sum_{l = 1}^N \widetilde  a_{mnl}\nabla \varphi_n(\x)  \cdot  \nabla v_l^{(p)}(\x)\Big|^2d\x
		+ \epsilon \|\varphi\|^2_{H^s(\Omega)^N}
		\label{3.5}
\end{multline}
for all $\varphi = [
		\begin{array}{lll}
			\varphi_1 &\dots &\varphi_N
		\end{array}
	]^{\rm T} \in H.$
Here, $\epsilon$ is a regularization parameter and $ \epsilon \|\varphi\|^2_{H^s(\Omega)^N}$ serves as a regularization term.
The ``Carleman parameters" $\lambda$ and $\beta$ and the regularization parameter $\epsilon$ will be chosen by a trial-error procedure, see Subsection \ref{sub42} for more details.
The following theorem holds.
\begin{Theorem}
 Let $H$ be defined as in \eqref{3.3}. Recursively define a sequence $\{{\bf v}^{(p)}\}_{p \geq 0}$ as in \eqref{3.4} where ${\bf v}^{(0)}$ is a vector-valued function in $H$.
	Assume the true solution $\bv^*$ to \eqref{2.14} is in $H$ with $\|\bv^*\|_{H^s(\Omega)^N} < M$.
	Fix $\beta > \beta_0$ where $\beta_0$ was defined in Lemma \ref{lemma carl}.
Let  $\lambda_0 = \lambda_0( \x_0,  \Omega, d, \beta)$ be the number in Lemma \ref{lemma carl}. Then, for all $\lambda \geq \lambda_0$,
	we have 
	\begin{multline}
		 \int_{\Omega}  e^{2\lambda r^{-\beta}(\x)}[|\nabla (\bv^{(p+1)} - \bv^*)|^2 + |\bv^{(p+1)} - \bv^*|^2 ]\,d\x
		+ \epsilon \|\bv^{(p+1)} - \bv^*\|_{H^s(\Omega)^N}^2
		\\
		\leq
		\eta^{p+1} \Big(\int_{\Omega}e^{2\lambda r^{-\beta}(\x)}	
			\big[|\nabla  \big(\bv^{(0)} - \bv^*\big)|^2
			+ |\bv^{(0)} - \bv^*|^2\big]d\x
			+ \epsilon \| \bv^{(0)} - \bv^*\|_{H^s(\Omega)^N}^2\Big)
			\\
			+  \frac{ \epsilon \eta}{1 - \eta}  \|\bv^*\|_{H^s(\Omega)^N}^2 
			\label{3.20}
	\end{multline}
    where $\eta \in (0, 1)$ is a positive constant depending only on $\x_0,$ $\Omega,$ $\beta$, $\{\Psi_n\}_{n = 1}^N$, $M$,  and $d$.
	\label{thm}
\end{Theorem}

\begin{Remark}
The conclusion of Theorem \ref{thm} is similar to those found in the convergence theorems outlined in \cite{LeCON2023, LeNguyen:jiip2022, Nguyen:AVM2023}. 
However, direct application of the convergence theorems in those publications to derive Theorem \ref{thm} is not possible. This limitation arises primarily because the nonlinearities present in the PDEs of \cite{LeCON2023, LeNguyen:jiip2022, Nguyen:AVM2023} must be Lipschitz, whereas the nonlinear function \[\bv \mapsto \sum_{n = 1}^N \sum_{l = 1}^N \widetilde  a_{mnl}\nabla v_n(\x)  \cdot  \nabla v_l(\x), \quad  m = 1, \dots, N,\] in \eqref{2.13} is not.
To overcome this obstacle, we substitute the Lipschitz condition with the boundedness assumption on $H$. 
\end{Remark}

\begin{proof}[Proof of Theorem \ref{thm}]

	Fix $p \geq 0$. Recall $J^{\bv^{(p)}}_{\lambda, \beta, \epsilon}$ as in \eqref{3.5}.
	Since ${\bf v}^{(p+1)}$ is the minimizer of $J^{\bv^{(p)}}_{\lambda, \beta, \epsilon}$ in $H$ and $H$ is a convex subset of $H^s(\Omega)^N$, for $0 < t < 1$, we have 
	\begin{equation} 
	\bv^{(p+1)} - t(   \bv^{(p+1)} - \bv^* )   
	= (1 - t) \bv^{(p+1)} + t \bv^*
	\in H.
	\end{equation} 
	Define $ {\bf h} = [
	\begin{array}{ccc}
		h_1& \dots & h_N
	\end{array}]^{\rm T} =\bv^{(p+1)} - \bv^*.$ We have
	\begin{equation}		
		\frac{J^{\bv^{(p)}}_{\lambda, \beta, \epsilon}\big(\bv^{(p+1)} - t {\bf h}\big) - J^{\bv^{(p)}}_{\lambda, \beta, \epsilon}\big(\bv^{(p+1)}\big)}{t} \geq 0.
		\label{3.11}
	\end{equation}
	Using \eqref{3.5} and \eqref{3.11}, we obtain
	\begin{multline}		
			\frac{1}{t}\sum_{m = 1}^N\Big[\Big\langle e^{2\lambda r^{-\beta}(\x)} \Big(2\Delta  v_m^{(p+1)}(\x) 
		+  2\sum_{n = 1}^N \nabla v_n^{(p+1)}(\x)  \cdot \widetilde  B_{mn} + 2\sum_{n = 1}^N \sum_{l = 1}^N \widetilde  a_{mnl}\nabla v_n^{(p+1)}(\x)  \cdot  \nabla v_l^{(p)}(\x)\Big)
		\\
		- t\Delta  h_m(\x) 
		-  t\sum_{n = 1}^N \nabla  h_n(\x)  \cdot \widetilde  B_{mn}- t\sum_{n = 1}^N \sum_{l = 1}^N \widetilde  a_{mnl}\nabla h_n(\x)  \cdot  \nabla v_l^{(p)}(\x),
		\\
		-t \Delta  h_m(\x) 
		-  t\sum_{n = 1}^N \nabla  h_n(\x)  \cdot \widetilde B_{mn} - t\sum_{n = 1}^N \sum_{l = 1}^N \widetilde a_{mnl}\nabla h_n(\x)  \cdot  \nabla v_l^{(p)}(\x)
		\Big\rangle_{L^2(\Omega)^N}
		\\
		+\epsilon \langle 2\bv^{(p+1)} - t{\bf h}, -t{\bf h}\rangle_{H^s(\Omega)^N}\Big]
		\geq 0.
		\label{3.14}
	\end{multline}
	Let $t$ in \eqref{3.14} tend to $0$. We have 
	\begin{multline}
		\sum_{m = 1}^N\Big\langle e^{2\lambda r^{-\beta}(\x)} \Big( \Delta  v_m^{(p+1)}(\x) 
		+  \sum_{n = 1}^N \nabla v_n^{(p+1)}(\x)  \cdot \widetilde B_{mn} + \sum_{n = 1}^N \sum_{l = 1}^N \widetilde a_{mnl}\nabla v_n^{(p+1)}(\x)  \cdot  \nabla v_l^{(p)}(\x)\Big),
		\\
	- \Delta  h_m(\x) 
		-  \sum_{n = 1}^N \nabla  h_n(\x)  \cdot \widetilde B_{mn} - \sum_{n = 1}^N \sum_{l = 1}^N \widetilde a_{mnl}\nabla h_n(\x)  \cdot  \nabla v_l^{(p)}(\x)
		\Big\rangle_{L^2(\Omega)^N}
		\\
		- \epsilon \langle \bv^{(p+1)}, {\bf h}\rangle_{H^s(\Omega)^N} \geq 0.
		\label{3.1616}
	\end{multline}
	or, equivalently,
	\begin{multline}
		\sum_{m = 1}^N\Big\langle e^{2\lambda r^{-\beta}(\x)} \Big( \Delta  v_m^{(p+1)}(\x) 
		+  \sum_{n = 1}^N \nabla v_n^{(p+1)}(\x)  \cdot \widetilde B_{mn} + \sum_{n = 1}^N \sum_{l = 1}^N \widetilde a_{mnl}\nabla v_n^{(p+1)}(\x)  \cdot  \nabla v_l^{(p)}(\x)\Big),
		\\
	+ \Delta  h_m(\x) 
		+  \sum_{n = 1}^N \nabla  h_n(\x)  \cdot \widetilde B_{mn} + \sum_{n = 1}^N \sum_{l = 1}^N \widetilde a_{mnl}\nabla h_n(\x)  \cdot  \nabla v_l^{(p)}(\x)
		\Big\rangle_{L^2(\Omega)^N}
		\\
		+ \epsilon \langle \bv^{(p+1)}, {\bf h}\rangle_{H^s(\Omega)^N} \leq 0.
		\label{3.16}
	\end{multline}
	On the other hand, since $\bv^*$ solves \eqref{2.13},
	\begin{multline}
		\sum_{m = 1}^N\Big\langle e^{2\lambda r^{-\beta}(\x)} \Big(\Delta  v_m^*(\x) 
		+  \sum_{n = 1}^N \nabla v_n^*(\x)  \cdot \widetilde B_{mn} + \sum_{n = 1}^N \sum_{l = 1}^N \widetilde a_{mnl}\nabla v_n^*(\x)  \cdot  \nabla v_l^*(\x)\Big),
		\\
	 \Delta  h_m(\x) 
		+  \sum_{n = 1}^N \nabla  h_n(\x)  \cdot \widetilde B_{mn} + \sum_{n = 1}^N \sum_{l = 1}^N \widetilde a_{mnl}\nabla h_n(\x)  \cdot  \nabla v_l^{(p)}(\x)
		\Big\rangle_{L^2(\Omega)^N}
		\\
		+ \epsilon \langle \bv^*, {\bf h}\rangle_{H^s(\Omega)^N} = \epsilon \langle \bv^*, {\bf h}\rangle_{H^s(\Omega)^N}.
		\label{3.17}
	\end{multline}
	Combining \eqref{3.16} and \eqref{3.17} gives
	\begin{multline}
		\sum_{m = 1}^N\Big\langle e^{2\lambda r^{-\beta}(\x)} \Big( \Delta  (v_m^{(p+1)}(\x) - v^*_m(\x))
		+  \sum_{n = 1}^N \nabla (v_n^{(p+1)}(\x) - v_n^*(\x))  \cdot \widetilde B_{mn} 
		+ \sum_{n = 1}^N \sum_{l = 1}^N \widetilde a_{mnl}\nabla v_n^{(p+1)}(\x)  \cdot  \nabla v_l^{(p)}(\x)
		\\
		-  \sum_{n = 1}^N \sum_{l = 1}^N \widetilde a_{mnl}\nabla v_n^*(\x)  \cdot  \nabla v_l^*(\x)\Big),
	 \Delta  h_m(\x) 
		+  \sum_{n = 1}^N \nabla  h_n(\x)  \cdot \widetilde B_{mn} 
		\\
		+ \sum_{n = 1}^N \sum_{l = 1}^N \widetilde a_{mnl}\nabla h_n(\x)  \cdot  \nabla v_l^{(p)}(\x)
		\Big\rangle_{L^2(\Omega)^N}
		+ \epsilon \langle \bv^{(p+1)} - \bv^*, {\bf h}\rangle_{H^s(\Omega)^N} 
		\leq -\epsilon \langle \bv^*, {\bf h}\rangle_{H^s(\Omega)^N}.
		\label{3.19}
	\end{multline}
	Recall that $\bh = \bv^{(p+1)} - \bv^*$. We obtain from \eqref{3.19} that
	\begin{multline*}
		\sum_{m = 1}^N \Big\langle e^{2\lambda r^{-\beta}(\x)} \Big( \Delta  h_m(\x)
		+  \sum_{n = 1}^N \nabla h_n(\x)  \cdot \widetilde B_{mn} 
		+ \sum_{n = 1}^N \sum_{l = 1}^N \widetilde a_{mnl}\nabla v_n^{(p+1)}(\x)  \cdot  \nabla[ v_l^{(p)}(\x) - v^*_l(\x)]
		\\
		+  \sum_{n = 1}^N \sum_{l = 1}^N \widetilde a_{mnl}\nabla h_n(\x)  \cdot  \nabla v_l^*(\x)\Big),
	 \Delta  h_m(\x) 
		+  \sum_{n = 1}^N \nabla  h_n(\x)  \cdot \widetilde B_{mn} 
		\\
		+ \sum_{n = 1}^N \sum_{l = 1}^N \widetilde a_{mnl}\nabla h_n(\x)  \cdot  \nabla v_l^{(p)}(\x)
		\Big\rangle_{L^2(\Omega)^N}
		+ \epsilon \| {\bf h}\|_{H^s(\Omega)^N}^2 
		\leq -\epsilon \langle \bv^*, {\bf h}\rangle_{H^s(\Omega)^N}.		
	\end{multline*}
	Thus, by using the inequality $2ab \leq \frac{1}{4} a^2 + 4b^2$, the H\"older inequality, and noting that $\|\bv^{(p+1)}\|_{H^s(\Omega)^N} \leq M$ and $\|\bv^*\|_{H^s(\Omega)^N} \leq M$ we get
	\begin{multline}
		\sum_{m = 1}^N\int_{\Omega}e^{2\lambda r^{-\beta}(\x)} \Big|
			 \Delta h_m(\x) 	+ \sum_{n = 1}^N \nabla h_n(\x)  \cdot \widetilde B_{mn} 
			\Big|
		^2d\x
		+\epsilon \| {\bf h}\|_{H^s(\Omega)^N}^2
		\\
		\leq 	C\Big(\int_{\Omega}e^{2\lambda r^{-\beta}(\x)}	
			|\nabla  \big(\bv^{(p)}(\x) - \bv^*(\x)\big)|^2d\x
			+ \int_{\Omega}e^{2\lambda r^{-\beta}(\x)}	
			|\nabla \bh|^2d\x
			\Big)
		+ \epsilon C \|\bv^*\|_{H^s(\Omega)^N}^2.
		\label{3.1414}
	\end{multline}
	Here, $C$ is a positive constant depending only on $\{\Psi_n\}_{n = 1}^N$ and $M$. 
	Using the inequality $(a + b)^2 \geq \frac{a^2}{2} - b^2$, we obtain from \eqref{3.1414} that
	\begin{multline}
		\int_{\Omega}e^{2\lambda r^{-\beta}(\x)} \Big|
			 \Delta \bh(\x)
			\Big|
		^2d\x
		+\epsilon \| {\bf h}\|_{H^s(\Omega)^N}^2
		\\
		\leq 	C\Big(\int_{\Omega}e^{2\lambda r^{-\beta}(\x)}	
			|\nabla  \big(\bv^{(p)}(\x) - \bv^*(\x)\big)|^2d\x
			+ \int_{\Omega}e^{2\lambda r^{-\beta}(\x)}	
			|\nabla \bh|^2d\x
			\Big)
		+ \epsilon C \|\bv^*\|_{H^s(\Omega)^N}^2.
		\label{3.15}
	\end{multline}
	From now on, we allow the universal constant $C$ to depend not only on $\{\Psi_n\}_{n = 1}^N$ and $M$ but also on $\x_0$, $\Omega,$ $\beta$, and $d$.
	Recall the number $\lambda_0$ as in Lemma \ref{lemma carl}.
	Applying the Carleman estimate \eqref{Car est} for each entry of $\bh$, we obtain
	\begin{align}
		\int_{\Omega} e^{2\lambda r^{-\beta}(\x)}\vert \Delta \bh \vert^2d\x
		&\geq
		 C\lambda  \int_{\Omega}  e^{2\lambda r^{-\beta}(\x)}\vert\nabla \bh(\x)\vert^2\,d\x
		+ C\lambda^3  \int_{\Omega}   e^{2\lambda r^{-\beta}(\x)}\vert \bh(\x)\vert^2\,d\x	\nonumber\\
		&\geq C\lambda  \int_{\Omega}  e^{2\lambda r^{-\beta}(\x)}[\vert\nabla \bh(\x)\vert^2  + |\bh(\x)|^2]\,d\x.
		\label{3.1616}	
	\end{align}
	It follows from \eqref{3.15} and \eqref{3.1616} that
	\begin{multline}
		\lambda \int_{\Omega}  e^{2\lambda r^{-\beta}(\x)}[\vert\nabla \bh(\x)\vert^2  + |\bh(\x)|^2]\,d\x
		+ \epsilon \| {\bf h}\|_{H^s(\Omega)^N}^2
		\\
		\leq 	C\Big(\int_{\Omega}e^{2\lambda r^{-\beta}(\x)}	
			|\nabla  \big(\bv^{(p)}(\x) - \bv^*(\x)\big)|^2d\x
			+ \int_{\Omega}e^{2\lambda r^{-\beta}(\x)}	
			|\nabla \bh|^2d\x
			\Big)
		+ \epsilon C \|\bv^*\|_{H^s(\Omega)^N}^2.
		\label{3.17}		
	\end{multline}
	Since \eqref{3.17} holds true for all $\lambda > \lambda_0$, we can choose $\lambda$ large such that the term $C \int_{\Omega} e^{2\lambda r^{-\beta}(\x)}|\nabla \bh(\x)|^2 d\x$ in the right-hand side of \eqref{3.17} is absorbed into the left-hand side. Recall that $\bh = \bv^{(p+1)} - \bv^*$. We obtain
	\begin{align}
	 \int_{\Omega} & e^{2\lambda r^{-\beta}(\x)}[|\nabla (\bv^{(p+1)} - \bv^*)|^2 + |\bv^{(p+1)} - \bv^*|^2 ]\,d\x
		+ \epsilon \|\bv^{(p+1)} - \bv^*\|_{H^s(\Omega)^N}^2 \nonumber
		\\
		&\leq 	\frac{C}{\lambda}\Big(\int_{\Omega}e^{2\lambda r^{-\beta}(\x)}	
			|\nabla  \big(\bv^{(p)} - \bv^*\big)|^2d\x
		+ \epsilon  \|\bv^*\|_{H^s(\Omega)^N}^2 \Big) \nonumber
		\\
		&\leq 	\frac{C}{\lambda}\Big(\int_{\Omega}e^{2\lambda r^{-\beta}(\x)}	
			\big[|\nabla  \big(\bv^{(p)} - \bv^*\big)|^2
			+ |\bv^{(p)} - \bv^*|^2\big]d\x
			+ \epsilon \| \bv^{(p)} - \bv^*\|_{H^s(\Omega)^N}^2
		+ \epsilon  \|\bv^*\|_{H^s(\Omega)^N}^2 \Big).
		\label{3.18}
	\end{align}
	Replacing $p$ by $p-1$ in \eqref{3.18}, we obtain
	\begin{align*}
		 \int_{\Omega}  e^{2\lambda r^{-\beta}(\x)}&[|\nabla (\bv^{(p+1)} - \bv^*)|^2 + |\bv^{(p+1)} - \bv^*|^2 ]\,d\x
		+ \epsilon \|\bv^{(p+1)} - \bv^*\|_{H^s(\Omega)^N}^2
		\\
		&\leq 
		\frac{C}{\lambda} \Big[
			\frac{C}{\lambda}\Big(\int_{\Omega}e^{2\lambda r^{-\beta}(\x)}	
			\big[|\nabla  \big(\bv^{(p-1)} - \bv^*\big)|^2
			+ |\bv^{(p-1)} - \bv^*|^2\big]d\x
			+ \epsilon \| \bv^{(p-1)} - \bv^*\|_{H^s(\Omega)^N}^2
			\\
		&\hspace{8cm}
			+ \epsilon  \|\bv^*\|_{H^s(\Omega)^N}^2 \Big)
			\Big)
		+ \epsilon  \|\bv^*\|_{H^s(\Omega)^N}^2 \Big)
		\Big]\\
		&= \Big(\frac{C}{\lambda}\Big)^2 \Big(\int_{\Omega}e^{2\lambda r^{-\beta}(\x)}	
			\big[|\nabla  \big(\bv^{(p-1)} - \bv^*\big)|^2
			+ |\bv^{(p-1)} - \bv^*|^2\big]d\x
			+ \epsilon \| \bv^{(p-1)} - \bv^*\|_{H^s(\Omega)^N}^2\Big)
			\\
		&\hspace{8cm}
			+ \Big(\Big(\frac{C}{\lambda}\Big)^2 + \frac{C}{\lambda}\Big)\epsilon  \|\bv^*\|_{H^s(\Omega)^N}^2.
	\end{align*}
	Continuing this process, we obtain
	\begin{multline}
		 \int_{\Omega}  e^{2\lambda r^{-\beta}(\x)}[|\nabla (\bv^{(p+1)} - \bv^*)|^2 + |\bv^{(p+1)} - \bv^*|^2 ]\,d\x
		+ \epsilon \|\bv^{(p+1)} - \bv^*\|_{H^s(\Omega)^N}^2
		\\
		\leq
		\Big(\frac{C}{\lambda}\Big)^{p+1} \Big(\int_{\Omega}e^{2\lambda r^{-\beta}(\x)}	
			\big[|\nabla  \big(\bv^{(0)} - \bv^*\big)|^2
			+ |\bv^{(0)} - \bv^*|^2\big]d\x
			+ \epsilon \| \bv^{(0)} - \bv^*\|_{H^s(\Omega)^N}^2\Big)
			\\
			+  \sum_{i = 1}^{p+1}\Big(\frac{C}{\lambda}\Big)^i \epsilon  \|\bv^*\|_{H^s(\Omega)^N}^2 .
			\label{3.19}
	\end{multline}
	Define $\eta = \frac{C}{\lambda} \in (0,1)$. The assertion \eqref{3.20} is directly deduced from \eqref{3.19}.	
\end{proof}

In theory, the estimate in \eqref{3.20} guarantee that when $p$ is sufficiently large, $\bv^{(p)}$ well-approximates the solution $\bv^*$ to \eqref{2.13}. However, in our numerical study, the approximation already meets our expectations when $p = 5$. 
As mentioned in Section \ref{sec2}, once \eqref{2.13} is solved, the computed solution to the inverse scattering problem follows. Our algorithm to solve this inverse scattering problem is presented in Algorithm \ref{alg}. Some details about implementation of Algorithm \ref{alg} will be presented in Subsection \ref{sub42}.
\begin{algorithm}[h!]
\caption{\label{alg}The Carleman contraction mapping method to compute the numerical solution to Problem \ref{p}}
	\begin{algorithmic}[1]
	\State \label{s_chooseN} Choose a cut-off number $N$.
	 Choose Carleman parameters $\x_0,$ $\beta$, and $\lambda$ and a regularization parameter $\epsilon$.
	 \State \label{s2} Choose an initial solution $\bv^{(0)}(\x) \in H$ and a number $P \in \mathbb N$. 
	 \For{$p = 0$ to $P-1$}
	 	\State \label{s4} Set $\bv^{(p+1)} = \underset{\varphi \in H}{\rm argmin} J^{\bv^{(p)}}_{\lambda, \beta, \epsilon}(\varphi)$ where $J^{\bv^{(p)}}_{\lambda, \beta, \epsilon}$ was defined in \eqref{3.4}.
	 \EndFor
	 \State   Set $\bv^{\rm comp} =  \bv^{(P)} =
	 \left[
	 	\begin{array}{cccc}
		v_1^{\rm comp} & v_2^{\rm comp} &\dots & v_N^{\rm comp}
		\end{array}
	 \right]^{\rm T}
	 $. By \eqref{2.7}, we compute $v^{\rm comp}(\x, \theta)$ using \begin{equation}
	v^{\rm comp}(\x, \theta) = \sum_{n = 1}^N v_n^{\rm comp}(\x) \Psi_n(\theta),
\end{equation}
for $\x \in \Omega,$ $\theta \in [0, 2\pi]$.
	\State 
	Due to \eqref{2.3}, a numerical solution to Problem \ref{p} can be computed via
	\begin{equation}
		c^{\rm comp}(\x) = \frac{1}{2\pi}\int_{0}^{2\pi}\Big|\Delta v(\x, \theta) + 2\ik \nabla v(\x, \theta) \cdot \hat \theta + k^2(\nabla v(\x, \theta))^2\Big|d\x d\theta  + 1
	\end{equation}
	for all $\x \in \Omega.$	 
\end{algorithmic}
\end{algorithm}

\begin{Remark}
For any given $\bv^{(p)} \in H$, it follows from \eqref{3.18} that $\bv^{(p+1)} = \underset{\varphi \in H}{\rm argmin} J^{\bv^{(p)}}_{\lambda, \beta, \epsilon}(\varphi)$ is closer to $\bv^*$, in comparison to $\bv^{(0)}$, with respect to the norm
\[
	\| \varphi \|_{\lambda, \beta, \epsilon}  =\Big[
	\int_{\Omega}  e^{2\lambda r^{-\beta}(\x)}[|\nabla \varphi|^2 + |\varphi|^2 ]\,d\x
		+ \epsilon \|\varphi\|_{H^s(\Omega)^N}^2 
		\Big]^{\frac{1}{2}}.
\]
So, the map $H \ni \bv \mapsto \underset{\varphi \in H}{\rm argmin} J^{\bv^{(p)}}_{\lambda, \beta, \epsilon}(\varphi)$ is a contractor toward $\bv^*$ with respect to the norm above. This contraction behavior of this map, together with the involvement of the Carleman weight function $e^{2\lambda r^{-\beta}(\x)}$, suggests the name of our method: the Carleman contraction mapping method.
\end{Remark}

\section{Numerical study}\label{sec4}

We display some numerical results based on Algorithm \ref{alg} in this section.
For simplicity, we only implement this algorithm in 2D.

\subsection{Data generation}

In this section, we discuss how to solve the forward problem to generate the data $f$ and $g$ for Problem \ref{p}.
Rather than solving \eqref{Hel_eqn} in the whole $\R^2$, which is complicated, we approximate this model by a Robin-type boundary problem for the Helmholtz equation
\begin{equation}
	\left\{
		\begin{array}{ll}
			\Delta u(\x, \theta) + k^2 c(\x) u(\x, \theta) = 0  &\x \in \Omega,
			\\
			\partial_{\nu} u_{\rm sc}(\x, \theta) - \ik  u_{\rm sc}(\x, \theta) = 0 &\x \in \partial \Omega,
		\end{array}
	\right.
	\label{4.1}
\end{equation}
for $\theta \in [0, 2\pi].$ Here, $u_{\rm sc}(\x, \theta) = u(\x, \theta) - u_{\rm inc}(\x, \theta)$ is the scattering wave.
This adjustment does not impact our methodology, as Algorithm \ref{alg} is constructed based on the form of the Helmholtz equation within $\Omega$ while the equation's form outside $\Omega$ is irrelevant.
In this paper, we set $\Omega = (-1, 1)^2$. On $\overline \Omega$, we arrange a uniform $64 \times 64$ grid of points
\[
	\mathcal G = \Big\{(x_i = -1 + (i - 1)d_\x, y_j = -1 + (j - 1)d_\x): 1 \leq i, j \leq 64 \Big\}
\]
where $d_\x = \frac{2}{63}$ is the step size in space.
We also discretize the interval of angle $[0, 2\pi]$ by the partition
\[
	\Theta = \big\{\theta_1 , \theta_2, \dots,  \theta_{N_\theta} \big\}
\] 
where $N_\theta = 150$, $d_\theta = 2\pi/149$, and $\theta_i = (i-1) d_\theta.$
It is not hard to verify that  for all $\theta \in [0, 2\pi],$ the scattering wave $u_{\rm sc}$ satisfies the following problem
\begin{equation}
	\left\{
		\begin{array}{ll}
			\Delta u_{\rm sc}(\x, \theta) + k^2 c(\x) u_{\rm sc}(\x, \theta) = -k^2(c(\x) - 1) u_{\rm inc}(\x, \theta)  &\x \in \Omega,
			\\
			\partial_{\nu} u_{\rm sc}(\x, \theta) - \ik  u_{\rm sc}(\x, \theta) = 0 &\x \in \partial \Omega.
		\end{array}
	\right.
	\label{4.2}
\end{equation}
In computation, we solve \eqref{4.2} by the finite difference method for $u_{\rm sc}$. Then, we compute $u$ by $u_{\rm sc} + u_{\rm inc}.$  
Then, we extract the ``noiseless" boundary data $f^*(\x, \theta)$ and $g^*(\x, \theta)$ for all $(\x, \theta) \in \mathcal G \times \Theta.$
The noisy data are given by
\[
	f = f^*(1 + \delta {\rm rand}), \quad
	g = g^*(1 + \delta {\rm rand})
\]
where rand is a function that generates random complex numbers with the range of both real part and imaginary part being $[-1, 1].$
The number $\delta$ is the noise level. In this paper, we set $\delta = 10\%.$

\subsection{The implementations}\label{sub42}

In this subsection, we discuss how we generate the simulated data and how we implement steps \ref{s_chooseN}, \ref{s2}, and \ref{s4} of Algorithm \ref{alg}. The implementation of other steps is straightforward.
	
\begin{figure}[h!]
	\centering
	\subfloat[Test 1]{\includegraphics[width=.22\textwidth]{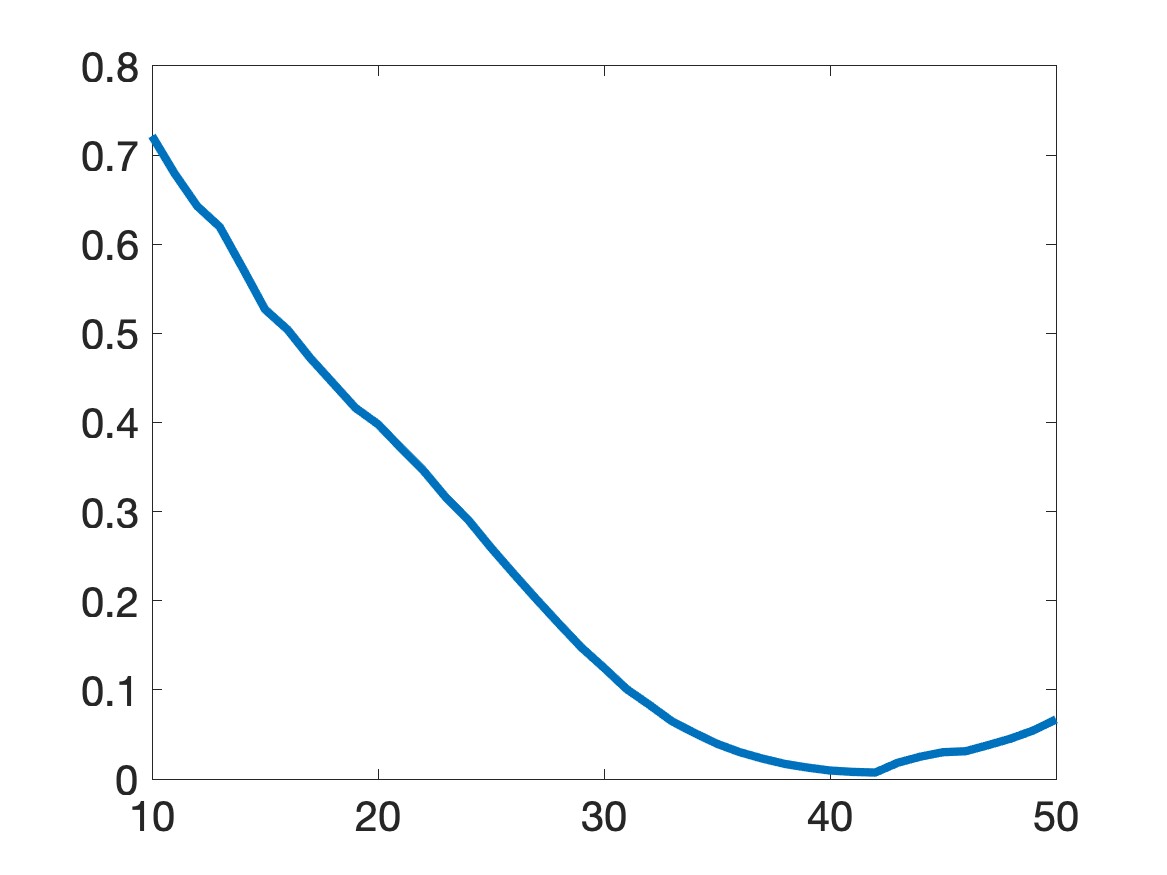}}
	\subfloat[Test 2]{\includegraphics[width=.22\textwidth]{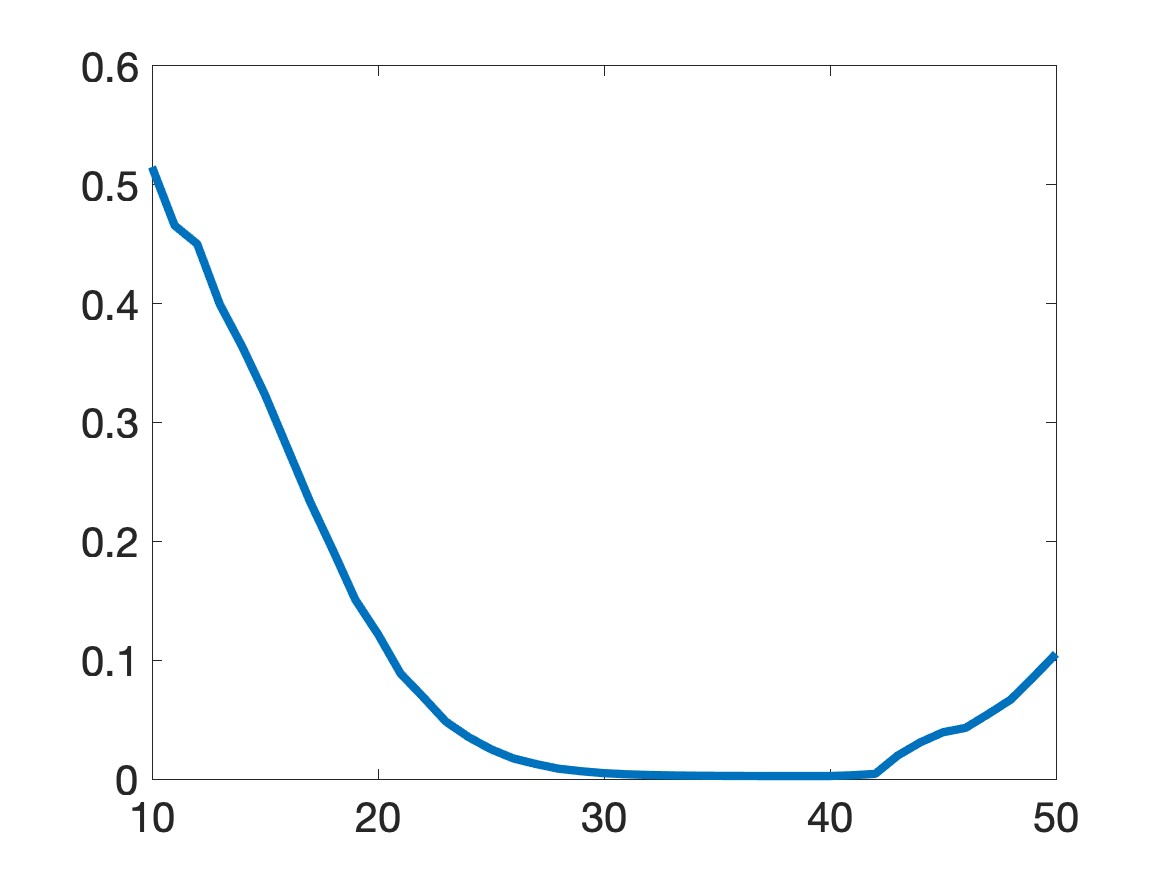}}
	\subfloat[Test 3]{\includegraphics[width=.22\textwidth]{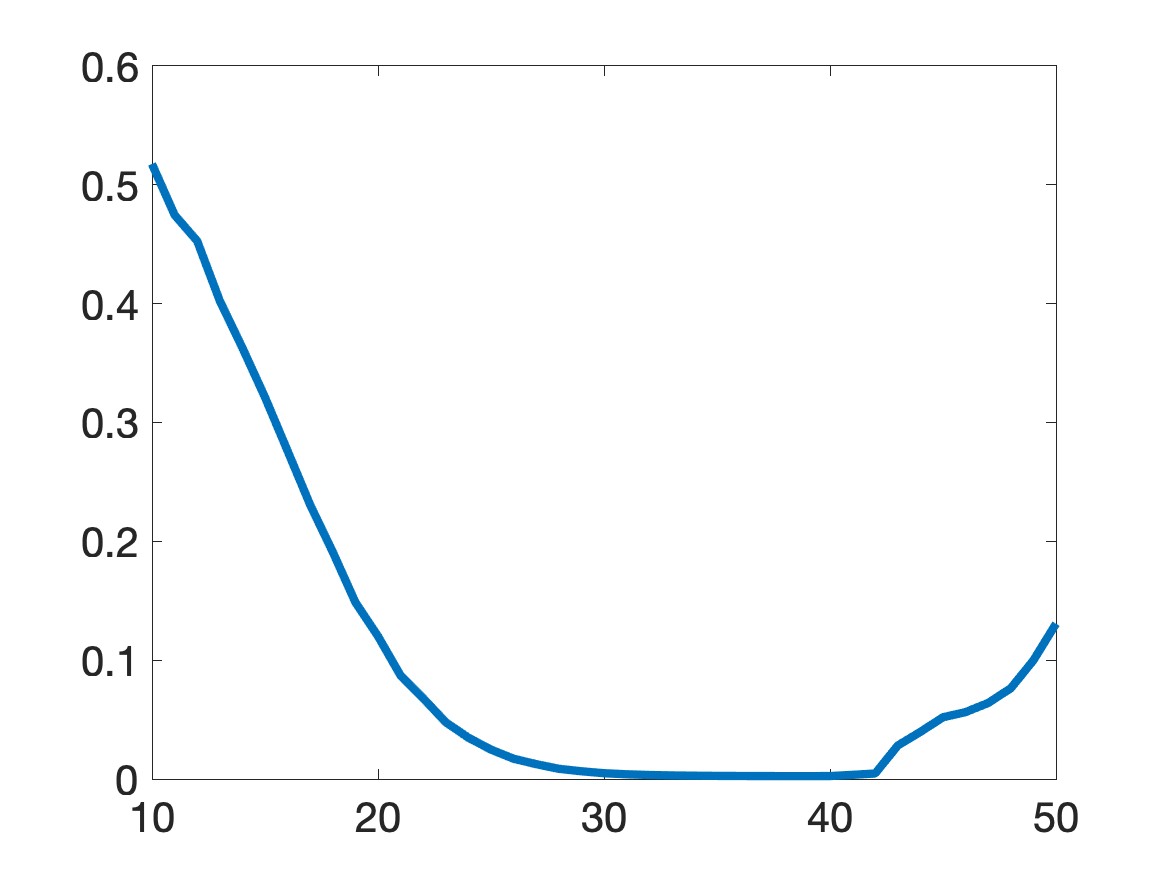}}
	\subfloat[Test 4]{\includegraphics[width=.22\textwidth]{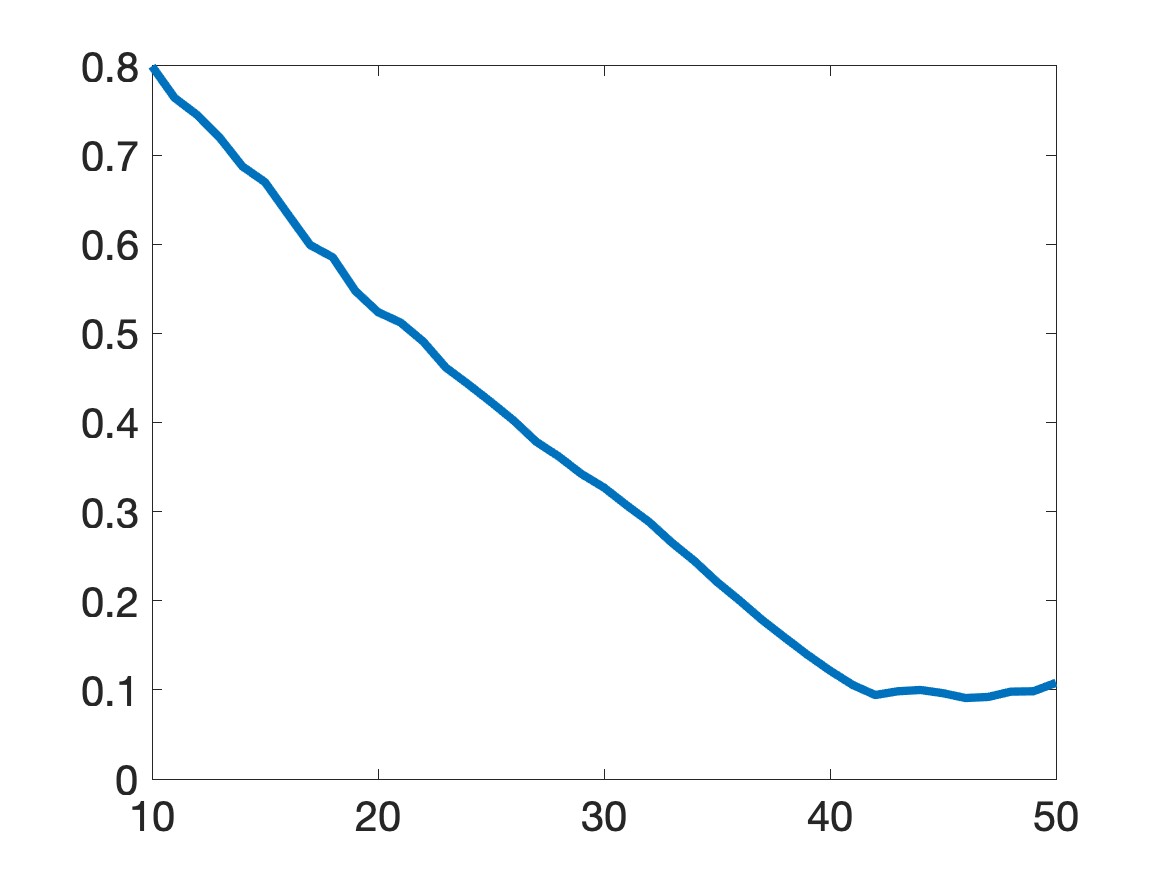}}
	\caption{\label{figChooseN} The graphs of the functions $e(N)$ are computed by the data, each containing 10\% noise, generated for all four numerical tests in Subsection \ref{sec_num_example}. These graphs suggest us to choose $N = 42$.	
}
\end{figure}

{\bf Step \ref{s_chooseN}: Choosing artificial parameters}. We use the information from the data $f$ to determine the appropriate cut-off parameter, $N$. Given that $\Omega$ is a square, it is convenient in numerical implementation to examine the data on one side, denoted as $\Gamma$, of its boundary $\partial \Omega$.
 For each $N$, we define the discrepancy between $f|_{\Gamma \times [0, 2\pi]}$ and its approximation obtained by truncating the Fourier series $\ds\sum_{n = 1}^N f_n(\mathbf{x}) \Psi_n(\theta)$, where $f_n(\x) = \ds \int_0^{2\pi} f(\x, \theta) \Psi_n(\theta)d\theta$. This difference is:
\[
	e(N) = \frac{\Big\|f - \ds\sum_{n = 1}^N f_n(\x) \Psi_n(\theta)\Big\|_{L^{2}(\Gamma \times [0, 2\pi])}}{\|f \|_{L^{2}(\Gamma \times [0, 2\pi])}}.
\]
Figure \ref{figChooseN} illustrates the plots of the function $e(N)$, where the data is sampled from the set \[\Gamma = \{(x, y = -1): |x| \leq 1\} \subset \partial \Omega.\] These data are generated for all numerical examples presented in Subsection \ref{sec_num_example}. In theory, the function $e(N)$ should decrease with respect to $N$. However, due to the high oscillation of $\Psi_n$ for large $n$, leading to errors in computing the integral in $\ds f_n(\x) = \int_0^{2\pi} f(\x, \theta) \Psi_n(\theta)d\theta$ using the finite difference method, this function experiences an increase at a certain critical value of $N$. That means the approximation is no longer true when $N$ is larger than this critical value. We select this critical value of $N$ as the cut-off parameter in Step \ref{s_chooseN}. Based on Figure \ref{figChooseN}, we choose $N = 42$.

A trial and error procedure chooses the other artificial parameters $\x_0, \beta, \lambda$, and $\epsilon$. We use Test 1 in Subsection \ref{sec_num_example} as a reference test, for which the true solution is known. These parameters are manually adjusted until the computed dielectric constant achieves an acceptable level of accuracy. Once determined, these parameters are then consistently applied across all subsequent tests. In our numerical study, we set $\mathbf{x}_0 = (0, -10)$, $\beta = 20$, $\lambda = 6$, and $\epsilon = 10^{-5.5}$.

{\bf Step \ref{s2}: Computing the initial solution}. Choosing the vector-valued function $\bv^{(0)} \in H$ is natural. To so so, we simply solve the over-determined boundary value problem
\begin{equation}
	\left\{
		\begin{array}{ll}
			\Delta v_m^{(0)} = 0 &\mbox{in } \Omega,\\
			v_m^{(0)}(\x) = F_m(\x) &\x \in \partial \Omega,\\
		\partial_{\nu} v_m^{(0)}(\x) = G_m(\x) &\x \in \partial \Omega
		\end{array}
	\right.
\end{equation}
by the regularized least square method, or so-called the quasi-reversibility method. That means we minimize the cost functional and set the minimizer as $\bv^{(0)}$
\[
	\bv^{(0)} = \underset{\varphi \in H}{{\rm argmin}} \int_{\Omega} |\Delta \varphi|^2 d\x.
\]
We find that this initial estimate for $\mathbf{v}^*$ may not be particularly close to the true solution. However, it is a robust starting point for the iterative process in Algorithm \ref{alg}.

Alternatively, a simpler choice for $\mathbf{v}^{(0)}$ is to set $\mathbf{v}^{(0)}$ as the zero vector everywhere. Despite the fact that $\mathbf{v}^{(0)} \equiv 0$ does not belong to $H$, we have observed that Algorithm \ref{alg} still yields surprisingly satisfactory numerical results.

{\bf Step \ref{s4}: Computing $\bv^{(p+1)}$ given the knowledge of $\bv^{(p)}$.} We observe that although $S$ is invertible, the matrix $S^{-1}$ contains ``singular" entries with large number when $N = 42$ as in our choice. This might cause some unnecessary errors in computation. We, therefore, rather than solve \eqref{2.13}, we solve its equivalent system in which the partial differential equations are replaced by \eqref{2.9}. The corresponding definition of $\bv^{(p+1)}$ is 
\begin{equation}
	\bv^{(p + 1)} =  \underset{\varphi \in H}{\rm argmin} \widetilde J^{\bv^{(p)}}_{\lambda, \beta, \epsilon}(\varphi)
\end{equation}
where 
\begin{multline}
	\widetilde J^{\bv^{(p)}}_{\lambda, \beta, \epsilon}(\varphi) =\sum_{m = 1}^N \int_{\Omega} e^{2\lambda r^{-\beta}(\x)}\Big|  \sum_{n = 1}^N s_{mn}\Delta  \varphi_n(\x) 
	\\
		+  \sum_{n = 1}^N \nabla  \varphi_n(\x)  \cdot  B_{mn} + \sum_{n = 1}^N \sum_{l = 1}^N   a_{mnl}\nabla \varphi_n(\x)  \cdot  \nabla v_l^{(p)}(\x)\Big|^2d\x
		+ \epsilon \|\varphi\|^2_{H^s(\Omega)^N}
		\label{4.5}
\end{multline}
for all $\varphi = [
		\begin{array}{lll}
			\varphi_1 &\dots &\varphi_N
		\end{array}
	]^{\rm T} \in H.$
	
	Note that the maps \[
		H \ni \varphi  \mapsto e^{\lambda r^{-\beta}(\x)} \sum_{n = 1}^N s_{mn}\Delta  \varphi_n(\x) 
	\\
		+  \sum_{n = 1}^N \nabla  \varphi_n(\x)  \cdot  B_{mn} + \sum_{n = 1}^N \sum_{l = 1}^N   a_{mnl}\nabla \varphi_n(\x)  \cdot  \nabla v_l^{(p)}(\x)
	\] 
	and
	\[
		H \ni \varphi \mapsto \sum_{|\alpha| \leq s} D^{\alpha} \varphi
	\]
	are linear.  Consequently, we can discretize all functions involved in \eqref{4.5} using a finite difference scheme. For a comprehensive understanding of the discretization process, we refer interested readers to \cite{Nguyen:CAMWA2020, Nguyens:jiip2020}.
The minimization procedure after discretion can be executed using the optimization packages available in Matlab. In our implementation, we employ the \texttt{lsqlin} command, which is a built-in tool in Matlab, for this purpose.

\subsection{Numerical examples}  \label{sec_num_example}

We display four numerical tests.

\begin{figure}[h!]
	\centering
	\subfloat[The function $c^{\rm true}$]{\includegraphics[width=.3\textwidth]{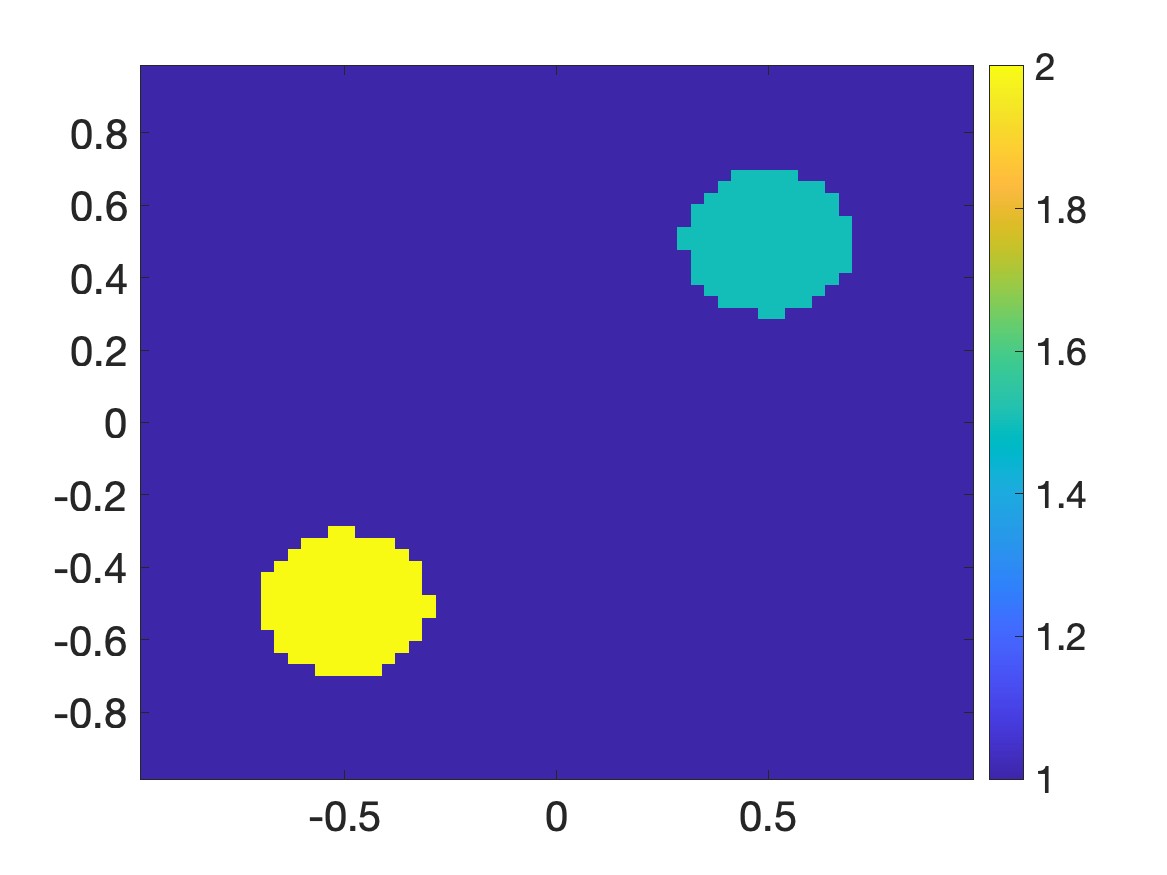}}
	\quad
	\subfloat[The function $c^{\rm comp}$]{\includegraphics[width=.3\textwidth]{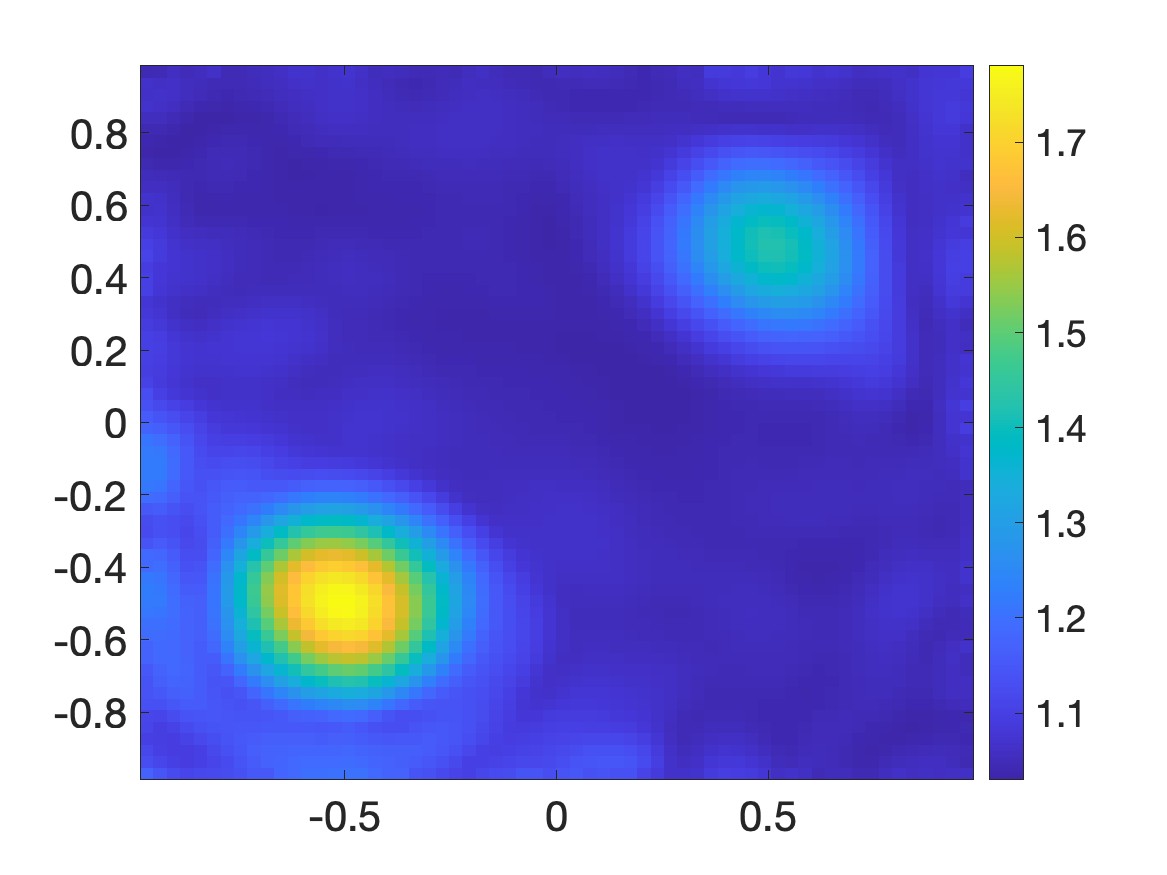}}
	\quad
	\subfloat[\label{fig1c} The consecutive difference $\frac{\|\bv^{(p)} - \bv^{(p-1)}\|_{L^2(\Omega)}}{\|\bv^{(p)}\|_{L^2(\Omega)}}$]{\includegraphics[width=.3\textwidth]{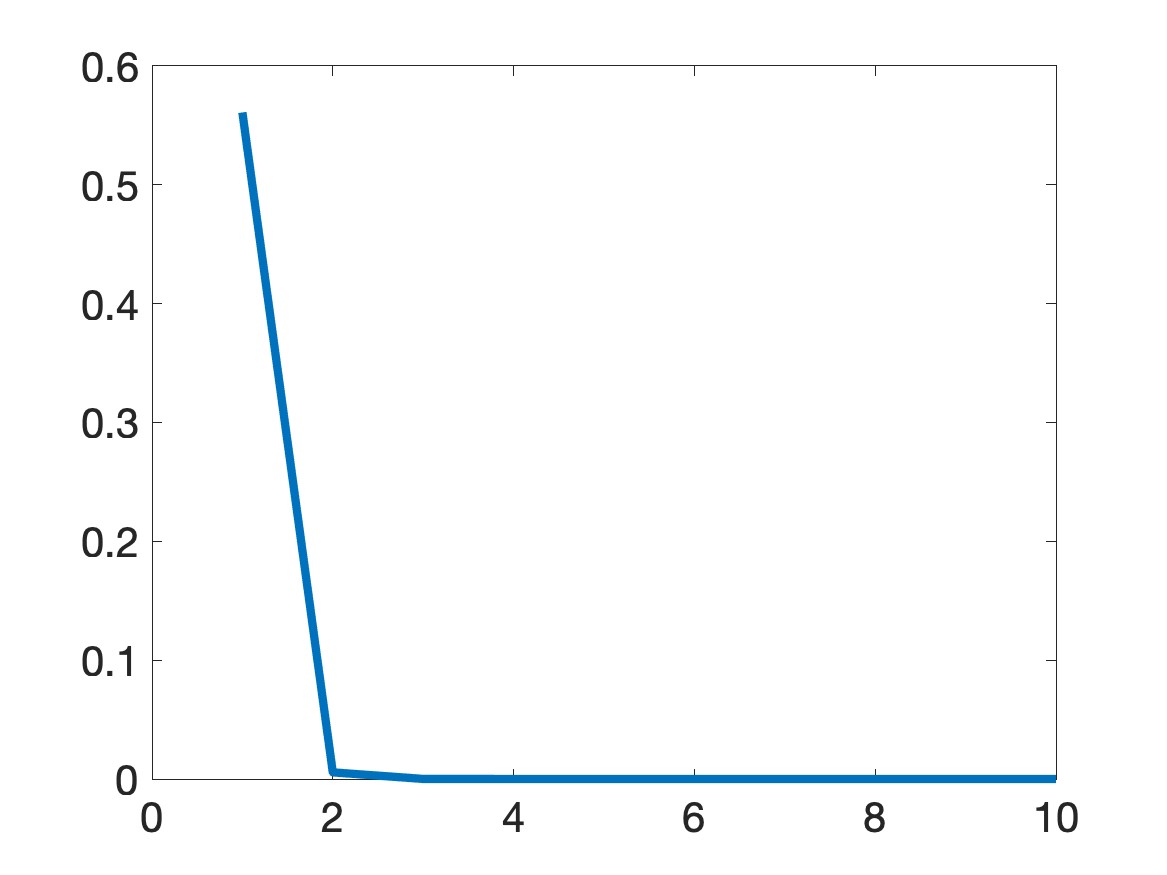}}
	\caption{\label{figtest1} The actual and computed spatial dielectric constants $c$ along with the $L^2$ relative difference of the computed $\mathbf{v}$ are presented. It is seen that the computed function $c$ meets the acceptable criteria. Additionally, Figure (c) provides numerical evidence supporting the rapid convergence of Algorithm \ref{alg}. Notably, it can be observed that only 3 iterations are adequate to achieve an accurate reconstruction. The data in this test is corrupted with 10\% of noise.}
\end{figure}

{\bf Test 1.} We test the case when ``two scatterers" with disk shape and with different values are presented inside $\Omega$. The function $c^{\rm true}$ is given by
\begin{equation}
	c^{\rm true}(x, y) =
	\left\{
		\begin{array}{ll}
			2 &(x + 0.5)^2 + (y + 0.5)^2 \leq 0.04,\\
			1.5 &(x - 0.5)^2 + (y - 0.5)^2 \leq 0.04,\\
			1 &\mbox{otherwise}. 
		\end{array}
	\right.
\end{equation}
In this test, the wave number $k$ is $3\pi.$
The true and constructed spatially distributed dielectric constants are displayed in Figure \ref{figtest1}.

In this test, the presence of two scatterers within $\Omega$ is clearly seen.  The errors in the computed results remain within acceptable bounds. The true value of $c$ within the lower scatterer is known to be 2, while the maximum value of $c$ calculated within the computed scatterer reaches 1.78, resulting in a relative error of 11.0\%. For the upper scatterer, the true $c$ value is 1.5, with the maximum computed $c$ value inside the computed scatterer being 1.42, yielding a relative error of 5.3\%. Both computed errors are compatible with the noise level $10\%$ present in the data. We highlight that Figure \ref{fig1c} serves as numerical verification of the convergence, proved in Theorem \ref{thm}.

\begin{figure}[h!]
	\centering
	\subfloat[The function $c^{\rm true}$]{\includegraphics[width=.3\textwidth]{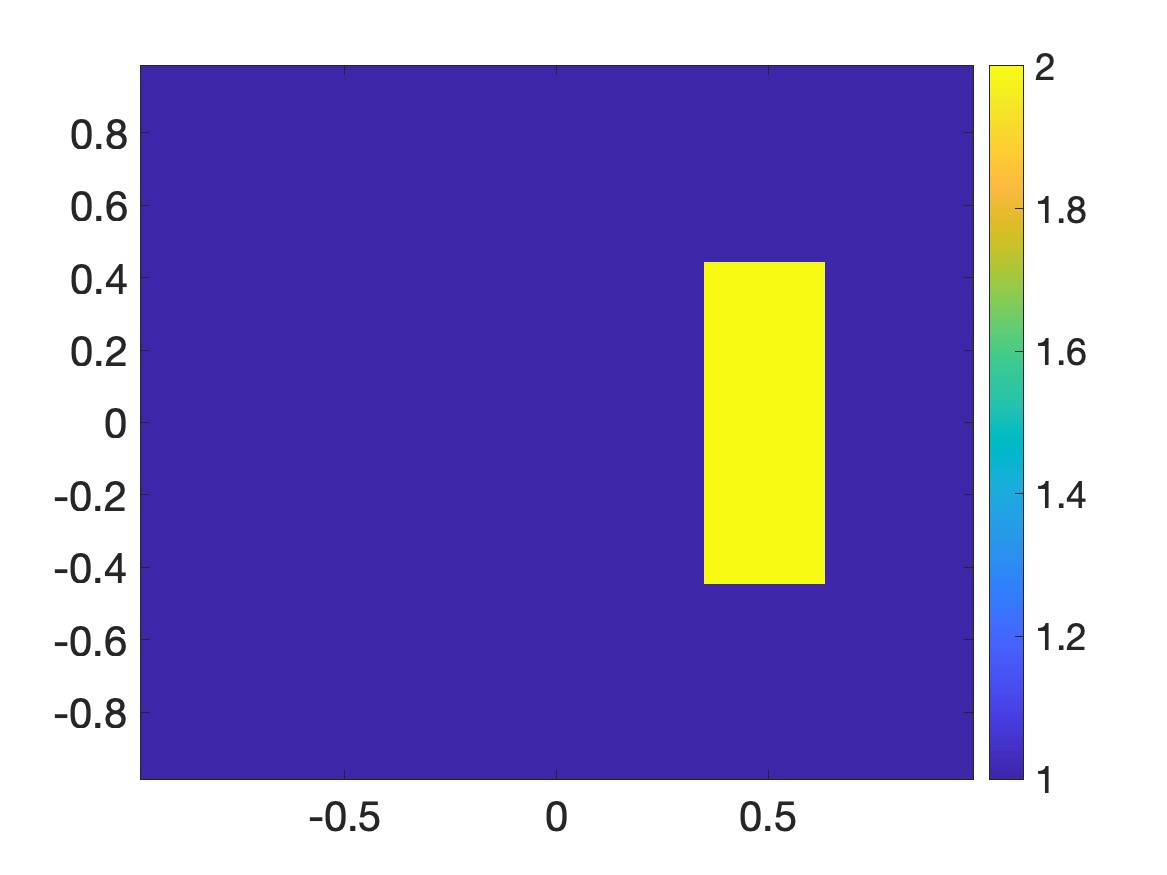}}
	\quad
	\subfloat[The function $c^{\rm comp}$]{\includegraphics[width=.3\textwidth]{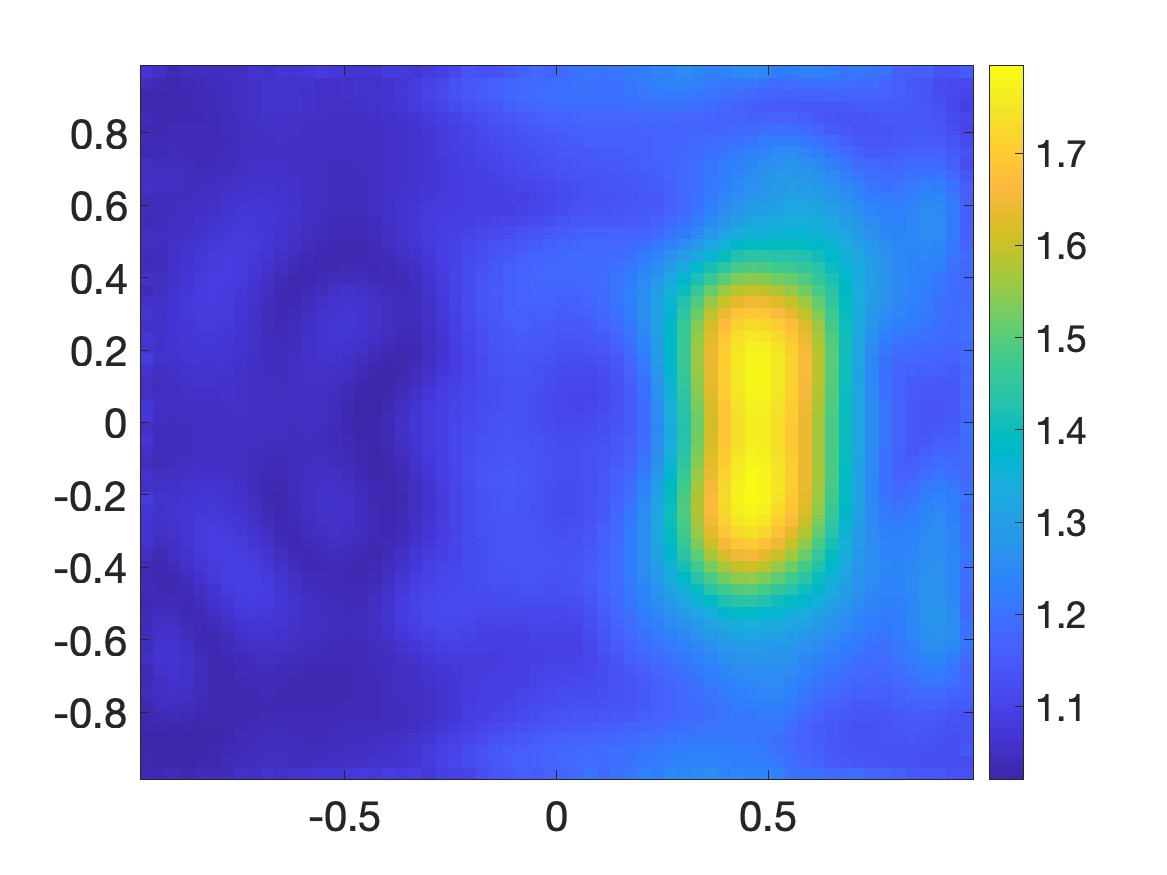}}
	\quad
	\subfloat[\label{fig2c} The consecutive difference $\frac{\|\bv^{(p)} - \bv^{(p-1)}\|_{L^2(\Omega)}}{\|\bv^{(p)}\|_{L^2(\Omega)}}$]{\includegraphics[width=.3\textwidth]{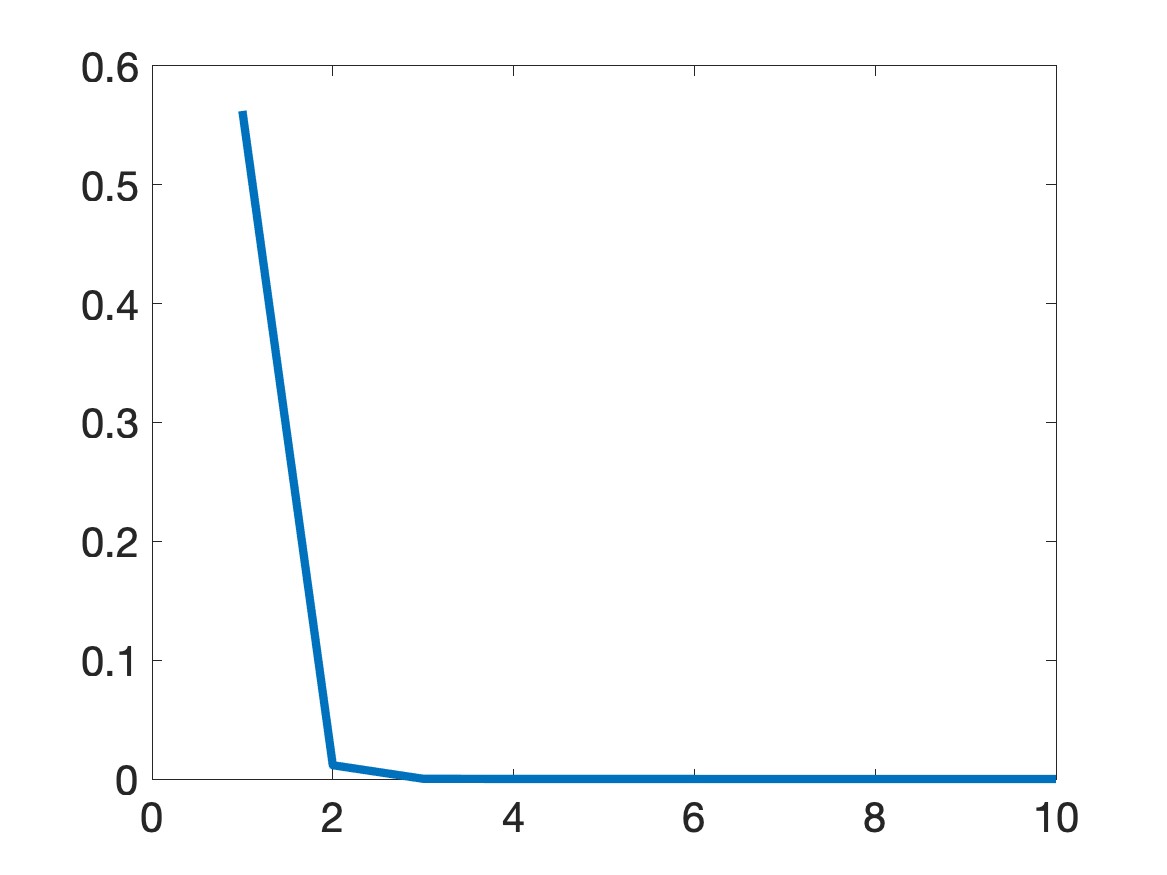}}
	\caption{\label{figtest2} The actual and computed spatial dielectric constants $c$ along with the $L^2$ relative difference of the computed $\mathbf{v}$ are presented. It is seen that the computed function $c$ meets the acceptable criteria. Additionally, Figure (c) provides numerical evidence supporting the rapid convergence of Algorithm \ref{alg}. Notably, it can be observed that only 3 iterations are adequate to achieve an accurate reconstruction. The data in this test is corrupted with 10\% of noise.}
\end{figure}

{\bf Test 2.} We test the case when ``one rectangle scatterer" is occluded inside $\Omega$. The function $c^{\rm true}$ is given by
\begin{equation}
	c^{\rm true}(x, y) =
	\left\{
		\begin{array}{ll}
			2 &\max\Big\{\frac{|x - 0.5|}{0.5}, \frac{|y|}{1.5}\Big\} \leq 0.3,\\
			1 &\mbox{otherwise}. 
		\end{array}
	\right.
\end{equation}
In this test, the wave number $k$ is $2\pi.$
The true and constructed spatially distributed dielectric constants are displayed in Figure \ref{figtest2}.

In this test,  the presence of the rectangular scatterer $\Omega$ is clearly seeable, with an acceptable reconstructed shape and location. Regarding the reconstructed value, the true value of $c$ in the scatterer is 2 while the maximum value of $c$ calculated within the reconstructed scatterer reaches 1.79, resulting in a relative error of 10.5\%.  Similarly to Test 1, Figure \ref{fig2c} not only visually confirms the convergence, as proven in Theorem \ref{thm}, but also serves as numerical validation of the convergence behavior detailed in the theorem.

\begin{figure}[h!]
	\centering
	\subfloat[The function $c^{\rm true}$]{\includegraphics[width=.3\textwidth]{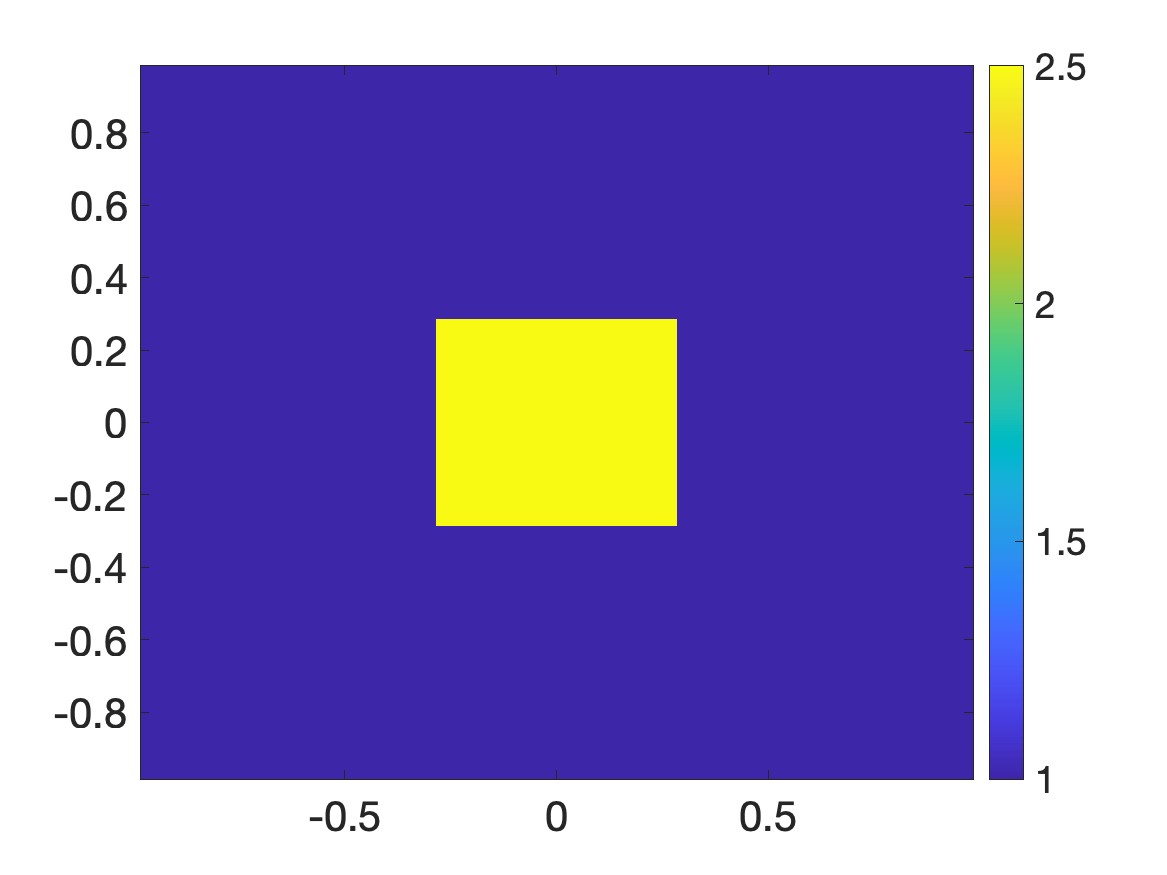}}
	\quad
	\subfloat[The function $c^{\rm comp}$]{\includegraphics[width=.3\textwidth]{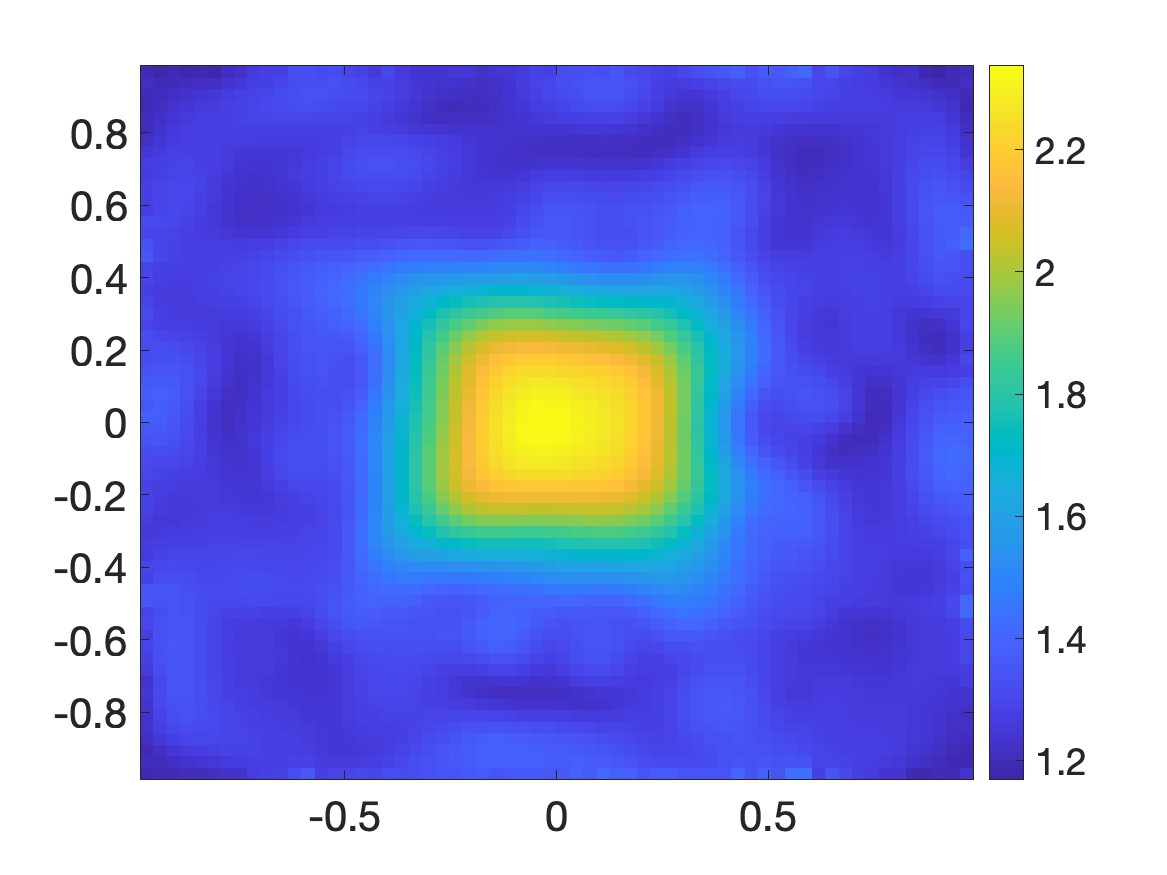}}
	\quad
	\subfloat[\label{fig3c} The consecutive difference $\frac{\|\bv^{(p)} - \bv^{(p-1)}\|_{L^2(\Omega)}}{\|\bv^{(p)}\|_{L^2(\Omega)}}$]{\includegraphics[width=.3\textwidth]{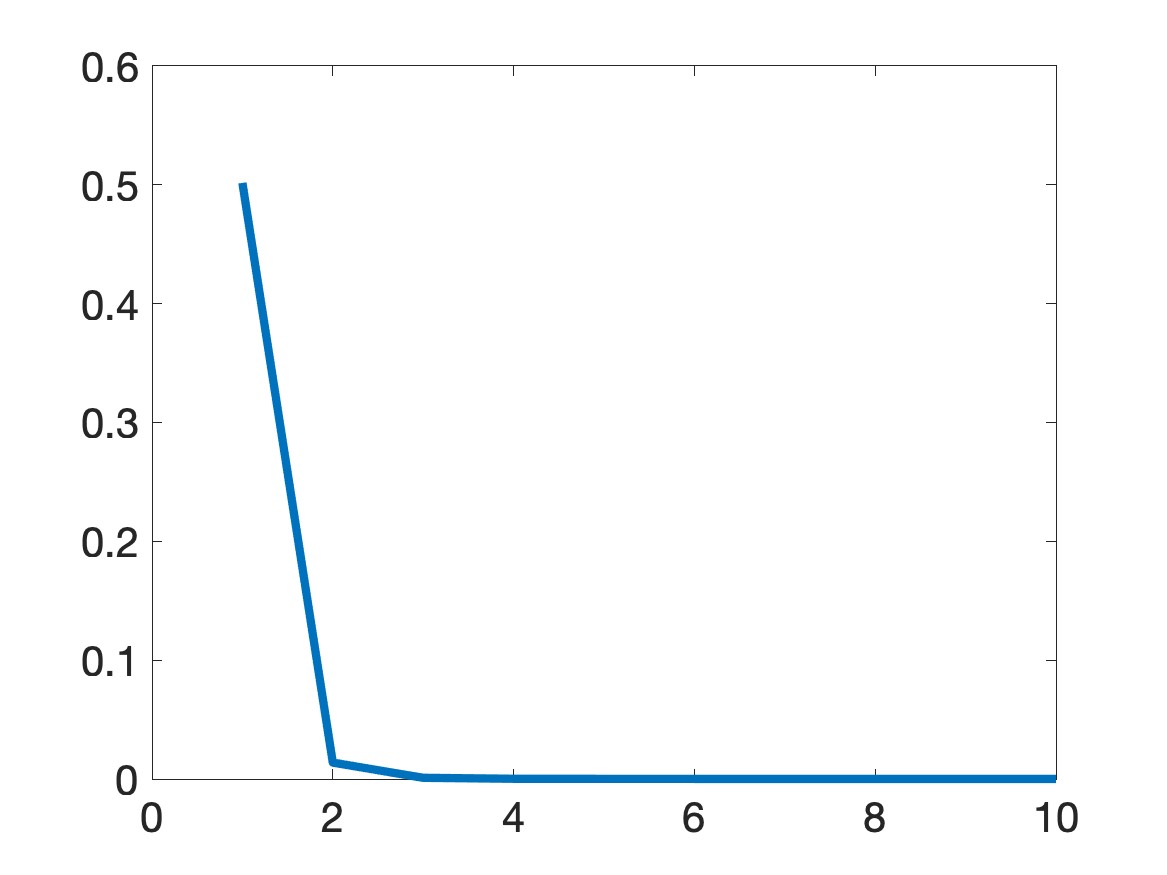}}
	\caption{\label{figtest3} The actual and computed spatial dielectric constants $c$ along with the $L^2$ relative difference of the computed $\mathbf{v}$ are presented. It is seen that the computed function $c$ meets the acceptable criteria. Additionally, Figure (c) provides numerical evidence supporting the rapid convergence of Algorithm \ref{alg}. Notably, it can be observed that only 3 iterations are adequate to achieve an accurate reconstruction. The data in this test is corrupted with 10\% of noise.}
\end{figure}

{\bf Test 3.} We test the case when ``one square scatterer" is placed at the center of $\Omega$, which is away from the measurement sites. The function $c^{\rm true}$ is given by
\begin{equation}
	c^{\rm true}(x, y) =
	\left\{
		\begin{array}{ll}
			2.5 &\max\{|x|, |y|\} \leq 0.09,\\
			1 &\mbox{otherwise}. 
		\end{array}
	\right.
\end{equation}
In this test, the wave number $k$ is $2\pi.$
The true and constructed spatially distributed dielectric constants are displayed in Figure \ref{figtest3}.

In this test,  the presence of the square scatterer $\Omega$ is clearly seeable, with an acceptable reconstructed shape and location at the center of $\Omega$. Regarding the reconstructed value, the true value of $c$ in the square scatterer is 2.5, while the maximum value of $c$ calculated within the reconstructed scatterer reaches 2.34. The relative error is 6.4\%.  Similarly to Test 1 and Test 2, Figure \ref{fig3c} plays the role of the evidence for the convergence behavior detailed in the theorem.

{\bf Test 4.}
We examine a test involving "fully occluded" scatterers. We consider two scatterers: one looking like a ring and the other resembling a disk placed within the first. The true function $c$ is defined as follows:
\begin{equation}
	c^{\rm true}(x, y) =
	\left\{
		\begin{array}{ll}
			1.5 & 0.5^2 < x^2 + y^2 < 0.7^2 \mbox{ or } x^2 + y^2 < 0.2^2\\
			1 &\mbox{otherwise}. 
		\end{array}
	\right.
\end{equation}

\begin{figure}[h!]
	\centering
	\subfloat[The function $c^{\rm true}$]{\includegraphics[width=.3\textwidth]{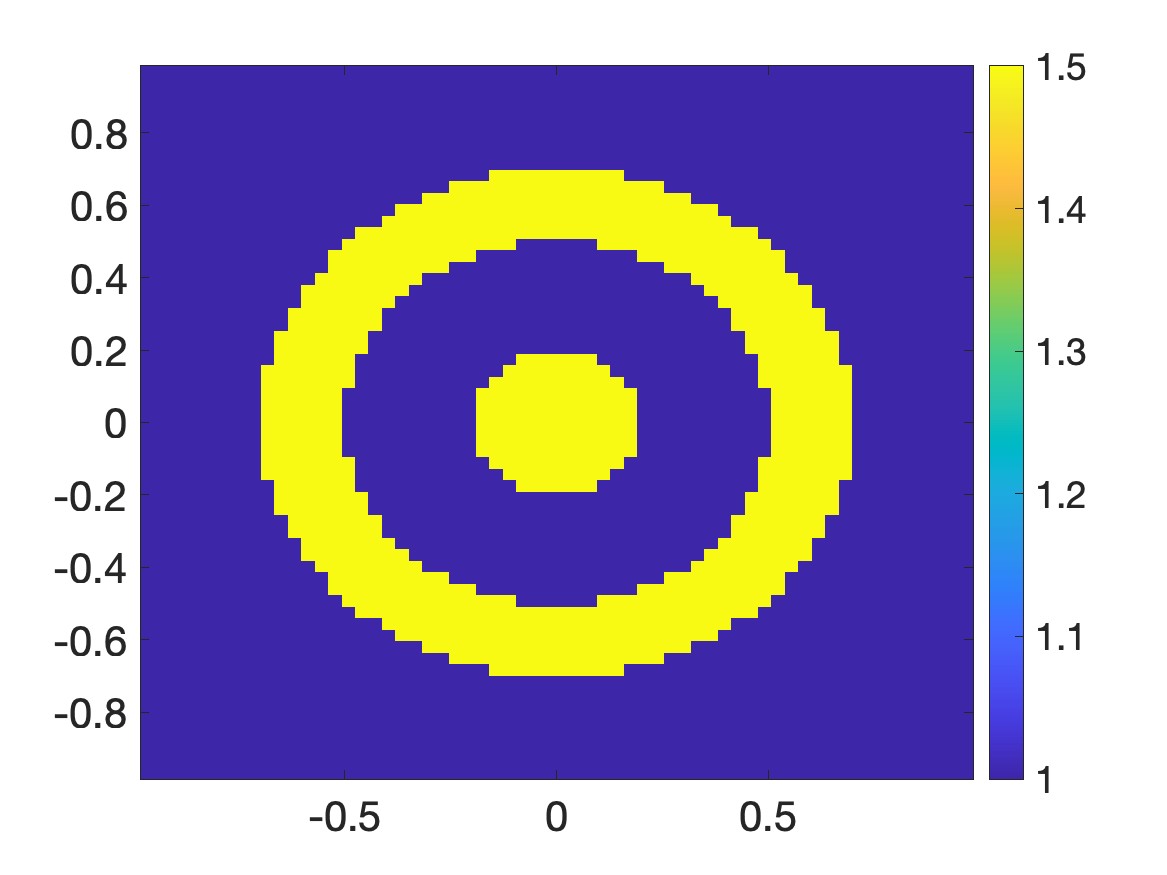}}
	\quad
	\subfloat[The function $c^{\rm comp}$]{\includegraphics[width=.3\textwidth]{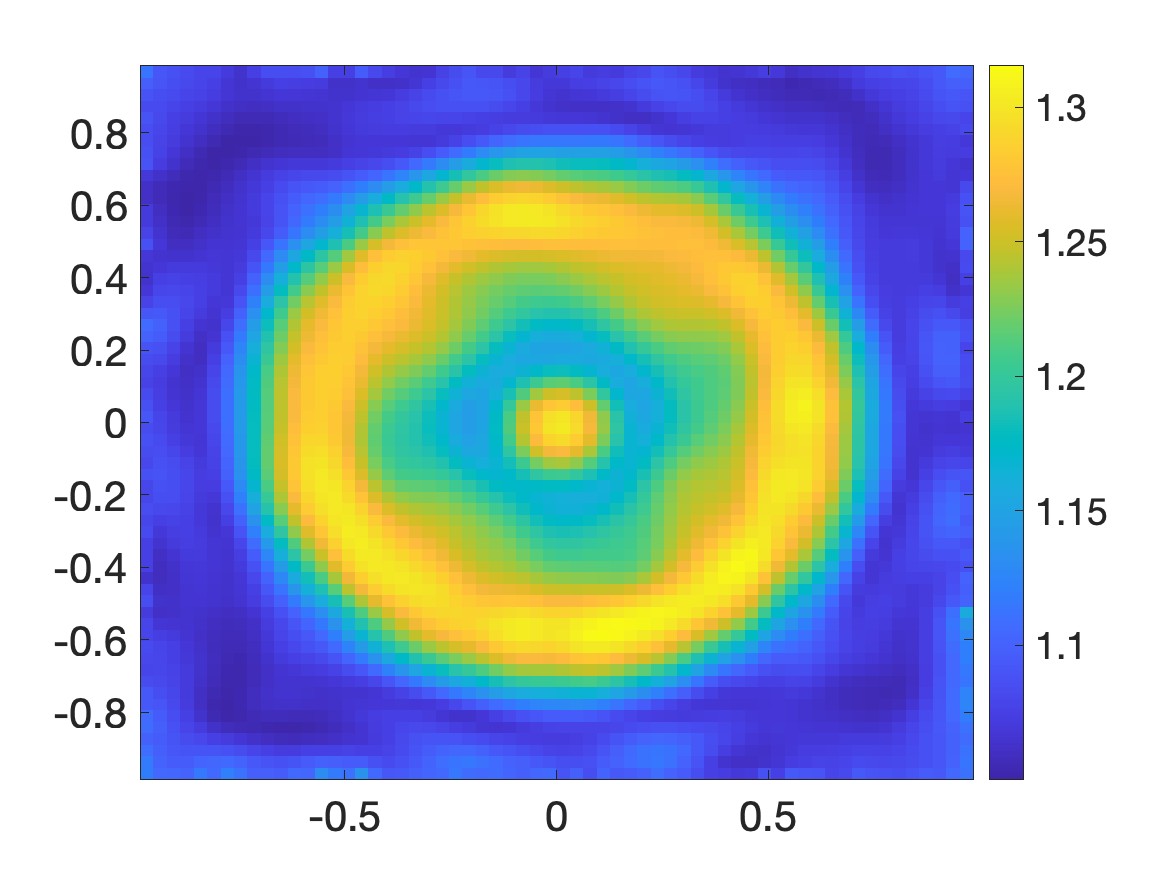}}
	\quad
	\subfloat[\label{fig4c} The consecutive difference $\frac{\|\bv^{(p)} - \bv^{(p-1)}\|_{L^2(\Omega)}}{\|\bv^{(p)}\|_{L^2(\Omega)}}$]{\includegraphics[width=.3\textwidth]{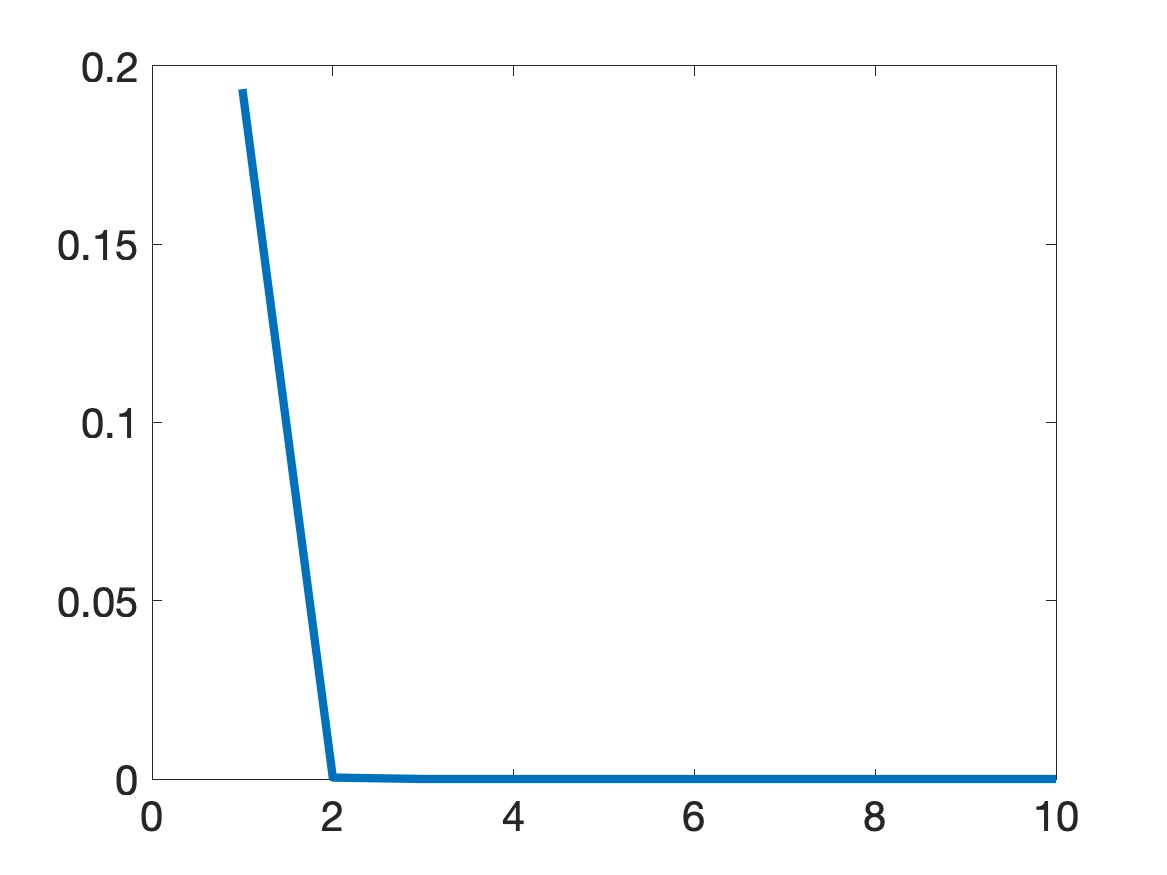}}
	\caption{\label{figtest4} The actual and computed spatial dielectric constants $c$ along with the $L^2$ relative difference of the computed $\mathbf{v}$ are presented. It is seen that the computed function $c$ meets the acceptable criteria. Additionally, Figure (c) provides numerical evidence supporting the rapid convergence of Algorithm \ref{alg}. Notably, it can be observed that only 3 iterations are adequate to achieve an accurate reconstruction. The data in this test is corrupted with 10\% of noise.}
\end{figure}

In this test, the wave number $k$ is $4\pi.$
The true and constructed spatially distributed dielectric constants are displayed in Figure \ref{figtest4}.

Observing both the outer and inner scatterers proves intriguing. Within the ring-shaped region, the true value of $c$ is known to be 1.5, while the maximum reconstructed value reaches 1.31, resulting in a relative error of 12.67\%. Similarly, within the inner scatterer, the true value of $c$ remains at 1.5, with the maximum reconstructed value reaching 1.30, yielding a relative error of 13.33\%. Notably, while this test presents a significant challenge due to the need to detect one inclusion within another, the relative errors remain consistent with the 10\% noise level. Once again, the convergence of our algorithm is verified both by Theorem \ref{thm} and by the observations presented in Figure \ref{fig4c}.

\section{Further discussions and concluding remarks}\label{sec5}

In the first part of this section, we delve into the computational cost of the Carleman contraction mapping method in Algorithm \ref{alg}. Remarkably, the computational expense of the  Carleman contraction mapping method employed in this study is notably modest. In particular, we only need to compute solutions for several linear systems to obtain the desired solutions to nonlinear problems in the tests above. As evident from Figures \ref{fig1c}--\ref{fig4c}, although solving 3 linear systems suffices, we execute this process 10 times to demonstrate convergence. Moreover, it is worth highlighting that the  Carleman contraction mapping method notably outperforms the Carleman convexification method introduced in our prior works such as \cite{VoKlibanovNguyen:IP2020, Khoaelal:IPSE2021, KhoaKlibanovLoc:SIAMImaging2020}. The computational Carleman convexification method incurs significant expense due to its involvement in the intensive procedure of minimizing a cost functional using the gradient descent method. On the other hand,  the  Carleman contraction mapping method does not require a precise initial guess of the true solution to the inverse scattering problem. Consequently, we assert that the  Carleman contraction mapping method is not only fast but also globally convergent.

The inverse scattering problem is widely recognized as severely ill-posed. Given the high noise levels typically encountered, one might question how we achieve satisfactory numerical results. The primary reason for our success lies in the truncation technique employed in \eqref{2.7}, which renders the original inverse scattering problem distinct from the approximation model \eqref{2.13}. An intriguing aspect emerges from this disparity: we effectively trade off equivalence for the challenge of solving a system of elliptic equations with Cauchy boundary data, a problem known for its stability. Consequently, our method demonstrates reduced sensitivity to noise, thereby facilitating robust numerical solutions.

In summary, this paper focuses on addressing the inverse scattering problem, which involves using waves generated through the interaction of plane waves with various directions and scatterers to detect these scatterers. These scatterers are characterized by their dielectric constant $c$.
Our approach to solving the inverse scattering problem comprises two key steps. Initially, we eliminate the unknown coefficient $c$ from the governing equation, resulting in a system of partial differential equations. Subsequently, we employ the  Carleman contraction mapping method to solve this system effectively. Once the solution to the aforementioned system is obtained, we can then compute the solution to the governing equation, thereby yielding the solution to the inverse scattering problem.

\section*{Declarations}

This manuscript is an original work. It has not been published previously and is not under consideration for publication elsewhere.
Both authors have contributed significantly and equally to the paper.

 \section*{Acknowledgement}
This work was partially supported by National Science Foundation grant DMS-2208159.

\bibliographystyle{plain}
\bibliography{../../../../mybib}

\begin{thebibliography}{10}

\bibitem{AbneyLeNguyenPeters}
R.~Abney, T.~T. Le, L.~H. Nguyen, and C.~Peters.
\newblock A {C}arleman-{P}icard approach for reconstructing zero-order
  coefficients in parabolic equations with limited data.
\newblock {\em preprint ArXiv:2309:14599}, 2023.

\bibitem{AmmariChowZou:sjap2016}
H.~Ammari, Y.~Chow, and J.~Zou.
\newblock Phased and phaseless domain reconstruction in inverse scattering
  problem via scattering coefficients.
\newblock {\em SIAM J. Appl. Math.}, 76:1000--1030, 2016.

\bibitem{AmmariKang:lnim2004}
H.~Ammari and H.~Kang.
\newblock {\em Reconstruction of Small Inhomogeneities from Boundary
  Measurements}, volume 1846.
\newblock Lecture Notes in Mathematics, Springer, 2004.

\bibitem{Bakushinsii:springer2004}
A.~B. Bakushinskii and M.~Y. Kokurin.
\newblock {\em Iterative Methods for Approximate Solutions of Inverse
  Problems}.
\newblock Springer, New York, 2004.

\bibitem{BaoLi:ip2005}
G.~Bao and P.~Li.
\newblock Inverse medium scattering for the {H}elmholtz equation at fixed
  frequency.
\newblock {\em Inverse Problems}, 21:1621--1641, 2005.

\bibitem{BaoLi:SIAM2005}
G.~Bao and P.~Li.
\newblock Inverse medium scattering problems for electromagnetic waves.
\newblock {\em SIAM J. Appl. Math.}, 65:2049--2066, 2005.

\bibitem{Bao:ip2015}
Gang Bao, P.~Li, J.~Lin, and F.~Triki.
\newblock Inverse scattering problems with multi-frequencies.
\newblock {\em Inverse Problems}, 31:093001, 2015.

\bibitem{BeilinaKlibanovBook}
L.~Beilina and M.~V. Klibanov.
\newblock {\em Approximate Global Convergence and Adaptivity for Coefficient
  Inverse Problems}.
\newblock Springer, New York, 2012.

\bibitem{Bleistein:ap1984}
N.~Bleistein.
\newblock {\em Mathematical Methods for Wave Phenomena}.
\newblock Academic Press, Orlando, 1984.

\bibitem{Tan1:ip2012}
T.~Bui-Thanh and O.~Ghattas.
\newblock Analysis of the {H}essian for inverse scattering problems: I.
  {I}nverse shape scattering of acoustic waves.
\newblock {\em Inverse Problems}, 28:055001, 2012.

\bibitem{Tan2:ip2012}
T.~Bui-Thanh and O.~Ghattas.
\newblock Analysis of the {H}essian for inverse scattering problems: {II.}
  {I}nverse medium scattering of acoustic waves.
\newblock {\em Inverse Problems}, 28:055002, 2012.

\bibitem{BukhgeimKlibanov:smd1981}
A.~L. Bukhgeim and M.~V. Klibanov.
\newblock Uniqueness in the large of a class of multidimensional inverse
  problems.
\newblock {\em Soviet Math. Doklady}, 17:244--247, 1981.

\bibitem{Burge2005}
M.~Burger and S.~Osher.
\newblock A survey on level set methods for inverse problems and optimal
  design.
\newblock {\em European J. of Appl. Math.}, 16:263--301, 2005.

\bibitem{Chavent:springer2009}
G.~Chavent.
\newblock {\em Nonlinear Least Squares for Inverse Problems: Theoretical
  Foundations and Step-by-Step Guide for Applications, Scientic Computation}.
\newblock Springer, New York, 2009.

\bibitem{Chen:ip1997}
Y.~Chen.
\newblock Inverse scattering via {H}eisenberg's uncertainty principle.
\newblock {\em Inverse Problems}, 13:253--282, 1997.

\bibitem{Chew:vnr1990}
W.~Chew.
\newblock {\em Waves and Fields in Inhomogeneous Media}.
\newblock Van Nostrand Reinhold, New York, 1990.

\bibitem{Colto1996}
D.~Colton and A.~Kirsch.
\newblock A simple method for solving inverse scattering problems in the
  resonance region.
\newblock {\em Inverse Problems}, 12:383--393, 1996.

\bibitem{ColtonKress:2013}
David Colton and Rainer Kress.
\newblock {\em Inverse acoustic and electromagnetic scattering theory.
  {A}pplied {M}athematical {S}ciences}.
\newblock Springer, New York, 3rd edition, 2013.

\bibitem{Devaney:cup2012}
A.~J. Devaney.
\newblock {\em Mathematical Foundations of Imaging, Tomography and Wavefield
  Inversion}.
\newblock Cambridge University Press, Cambridge, 2012.

\bibitem{Engl:Kluwer1996}
H.~W. Engl, M.~Hanke, and A.~Neubauer.
\newblock {\em Regularization of Inverse Problems, Mathematics and its
  Applications}.
\newblock Kluwer Academic Publishers Group, Dordrecht, 1996.

\bibitem{Gonch2013}
A.~V. Goncharsky and S.~Y. Romanov.
\newblock Supercomputer technologies in inverse problems of ultrasound
  tomography.
\newblock {\em Inverse Problems}, 29:075004, 2013.

\bibitem{HarrisLiem:SIAM2020}
I.~Harris and D-L. Nguyen.
\newblock Orthogonality sampling method for the electromagnetic inverse
  scattering problem.
\newblock {\em SIAM Journal on Scientific Computing}, 42:B722--B737, 2020.

\bibitem{VoKlibanovNguyen:IP2020}
V.~A. Khoa, G.~W. Bidney, M.~V. Klibanov, L.~H. Nguyen, L.~Nguyen, A.~Sullivan,
  and V.~N. Astratov.
\newblock Convexification and experimental data for a {3D} inverse scattering
  problem with the moving point source.
\newblock {\em Inverse Problems}, 36:085007, 2020.

\bibitem{Khoaelal:IPSE2021}
V.~A. Khoa, G.~W. Bidney, M.~V. Klibanov, L.~H. Nguyen, L.~Nguyen, A.~Sullivan,
  and V.~N. Astratov.
\newblock An inverse problem of a simultaneous reconstruction of the dielectric
  constant and conductivity from experimental backscattering data.
\newblock {\em Inverse Problems in Science and Engineering}, 29(5):712--735,
  2021.

\bibitem{KhoaKlibanovLoc:SIAMImaging2020}
V.~A. Khoa, M.~V. Klibanov, and L.~H. Nguyen.
\newblock Convexification for a 3{D} inverse scattering problem with the moving
  point source.
\newblock {\em SIAM J. Imaging Sci.}, 13(2):871--904, 2020.

\bibitem{Kirsc1998}
A.~Kirsch.
\newblock Characterization of the shape of a scattering obstacle using the
  spectral data of the far field operator.
\newblock {\em Inverse Problems}, 14:1489--1512, 1998.

\bibitem{Kirsch:aa2017}
A.~Kirsch.
\newblock Remarks on the {B}orn approximation and the factorization method.
\newblock {\em Appl. Anal.}, 96:70--84, 2017.

\bibitem{Klibanov:jiip2017}
M.~V. Klibanov.
\newblock Convexification of restricted {D}irichlet to {N}eumann map.
\newblock {\em J. Inverse and Ill-Posed Problems}, 25(5):669--685, 2017.

\bibitem{KlibanovIoussoupova:SMA1995}
M.~V. Klibanov and O.~V. Ioussoupova.
\newblock Uniform strict convexity of a cost functional for three-dimensional
  inverse scattering problem.
\newblock {\em SIAM J. Math. Anal.}, 26:147--179, 1995.

\bibitem{KlibanovLeNguyenIPI2022}
M.~V. Klibanov, T.~T. Le, L.~H. Nguyen, A.~Sullivan, and L.~Nguyen.
\newblock Convexification-based globally convergent numerical method for a 1{D}
  coefficient inverse problem with experimental data.
\newblock {\em Inverse Problems and Imaging}, 16:1579--1618, 2022.

\bibitem{KlibanovLiBook}
M.~V. Klibanov and J.~Li.
\newblock {\em Inverse Problems and Carleman Estimates: Global Uniqueness,
  Global Convergence and Experimental Data}.
\newblock De Gruyter, 2021.

\bibitem{KlibanovLiemLoc:SIAMjap2019}
M.~V. Klibanov, D-L Nguyen, and L.~H. Nguyen.
\newblock A coefficient inverse problem with a single measurement of phaseless
  scattering data.
\newblock {\em SIAM Journal of Applied Mathematics}, 79:1--27, 2019.

\bibitem{KlibanovNguyenTran:JCP2022}
M.~V. Klibanov, L.~H. Nguyen, and H.~V. Tran.
\newblock Numerical viscosity solutions to {H}amilton-{J}acobi equations via a
  {C}arleman estimate and the convexification method.
\newblock {\em Journal of Computational Physics}, 451:110828, 2022.

\bibitem{Klibanov:jiip2023}
M.~V. Klibanov and V.~G. Romanov.
\newblock A h\"older stability estimate for a 3d coefficient inverse problem
  for a hyperbolic equation with a plane wave.
\newblock {\em Journal of Inverse and Ill-posed Problems}, 31:223--242, 2023.

\bibitem{Langenberg:1987}
K.~J. Langenberg.
\newblock Applied inverse problems for acoustic, electromagnetic and elastic
  wave scattering.
\newblock In Sabatier, editor, {\em Basic Methods of Tomography and Inverse
  Problems}, pages 127--467. Adam Hilger, 1987.

\bibitem{Lavrentiev:AMS1986}
M.~M. Lavrent'ev, V.~G. Romanov, and S.~P. Shishat$\cdot$ski\u{i}.
\newblock {\em Ill-Posed Problems of Mathematical Physics and Analysis}.
\newblock Translations of Mathematical Monographs. AMS, Providence: RI, 1986.

\bibitem{LeCON2023}
T.~T. Le.
\newblock Global reconstruction of initial conditions of nonlinear parabolic
  equations via the {C}arleman-contraction method.
\newblock In D-L. Nguyen, L.~H. Nguyen, and T-P. Nguyen, editors, {\em Advances
  in Inverse problems for Partial Differential Equations}, volume 784 of {\em
  Contemporary Mathematics}, pages 23--42. American Mathematical Society, 2023.

\bibitem{LeNguyen:jiip2022}
T.~T. Le and L.~H. Nguyen.
\newblock A convergent numerical method to recover the initial condition of
  nonlinear parabolic equations from lateral {C}auchy data.
\newblock {\em Journal of Inverse and Ill-posed Problems}, 30(2):265--286,
  2022.

\bibitem{LeNguyenTran:CAMWA2022}
T.~T. Le, L.~H. Nguyen, and H.~V. Tran.
\newblock A {C}arleman-based numerical method for quasilinear elliptic
  equations with over-determined boundary data and applications.
\newblock {\em Computers and Mathematics with Applications}, 125:13--24, 2022.

\bibitem{LeNguyenNguyenPark}
T.~T. Le, L.~V. Nguyen, L.~H. Nguyen, and H.~Park.
\newblock The time dimensional reduction method to determine the initial
  conditions without the knowledge of damping coefficients.
\newblock {\em preprint arXiv:2308.13152}, 2023.

\bibitem{LiLiuWang:jcp2014}
J.~Li, H.~Liu, and Q.~Wang.
\newblock Enhanced multilevel linear sampling methods for inverse scattering
  problems.
\newblock {\em J. Comput. Phys.}, 257:554--571, 2014.

\bibitem{LiLiuZou:smms2014}
Z.~Li, H.~Liu, and J.~Zou.
\newblock Locating multiple multiscale acoustic scatterers.
\newblock {\em SIAM Multiscale Model. Simul.}, 12:927--952, 2014.

\bibitem{Moskow:Ip2008}
S.~Moskow and J.~Schotland.
\newblock Convergence and stability of the inverse {B}orn series for diffuse
  waves.
\newblock {\em Inverse Problems}, 24:065004, 2008.

\bibitem{NguyenNguyenTruong:camwa2022}
D-L Nguyen, L.~H. Nguyen, and T.~Truong.
\newblock The {C}arleman-based contraction principle to reconstruct the
  potential of nonlinear hyperbolic equations.
\newblock {\em Computers and Mathematics with Applications}, 128:239--248,
  2022.

\bibitem{Liem:jiip2022}
D-L. Nguyen and T.~Truong.
\newblock Imaging of bi-anisotropic periodic structures from electromagnetic
  near field data.
\newblock {\em Journal of Inverse and Ill-posed Problems}, 30:205--219, 2022.

\bibitem{MinhLoc:tams2015}
H.~M. Nguyen and L.~H. Nguyen.
\newblock Cloaking using complementary media for the {H}elmholtz equation and a
  three spheres inequality for second order elliptic equations.
\newblock {\em Transaction of the American Mathematical Society}, 2:93--112,
  2015.

\bibitem{LocNguyen:ip2019}
L.~H. Nguyen.
\newblock An inverse space-dependent source problem for hyperbolic equations
  and the {L}ipschitz-like convergence of the quasi-reversibility method.
\newblock {\em Inverse Problems}, 35:035007, 2019.

\bibitem{Nguyen:CAMWA2020}
L.~H. Nguyen.
\newblock A new algorithm to determine the creation or depletion term of
  parabolic equations from boundary measurements.
\newblock {\em Computers and Mathematics with Applications}, 80:2135--2149,
  2020.

\bibitem{Nguyen:AVM2023}
L.~H. Nguyen.
\newblock The {C}arleman contraction mapping method for quasilinear elliptic
  equations with over-determined boundary data.
\newblock {\em Acta Mathematica Vietnamica}, 48:401--422, 2023.

\bibitem{NguyenLiKlibanov:2019}
L.~H. Nguyen, Q.~Li, and M.~V. Klibanov.
\newblock A convergent numerical method for a multi-frequency inverse source
  problem in inhomogenous media.
\newblock {\em Inverse Problems and Imaging}, 13:1067--1094, 2019.

\bibitem{Nguyens:jiip2020}
P.~M. Nguyen and L.~H. Nguyen.
\newblock A numerical method for an inverse source problem for parabolic
  equations and its application to a coefficient inverse problem.
\newblock {\em Journal of Inverse and Ill-posed Problems}, 38:232--339, 2020.

\bibitem{Soumekh:SAR}
M.~Soumekh.
\newblock {\em Synthetic {A}perture {R}adar Signal Processing with {MATLAB
  A}lgorithms
  (https://www.mathworks.com/matlabcentral/fileexchange/2188-synthetic-aperture-radar-signal-processing-with-matlab-algorithms),
  {MATLAB Central File Exchange}.}
\newblock John Wiley \& Sons, 1999.

\bibitem{Tihkonov:kapg1995}
A.~N. Tikhonov, A.~Goncharsky, V.~V. Stepanov, and A.~G. Yagola.
\newblock {\em Numerical Methods for the Solution of Ill-Posed Problems}.
\newblock Kluwer Academic Publishers Group, Dordrecht, 1995.

\end{thebibliography}

\end{document}